\documentclass[a4paper, 12 pt]{article}

\usepackage[T2A]{fontenc}
\usepackage[cp1251]{inputenc}
\usepackage[english]{babel}
\usepackage[tbtags]{amsmath}
\usepackage{amsfonts,amssymb}

\usepackage{mathrsfs}

\usepackage{geometry} % Меняем поля страницы
\begin{document}
\title{Generalized Riordan groups and zero generalized Pascal matrices}
 \author{E. Burlachenko}
 \date{}

 \maketitle
\begin{abstract}
The generalized Riordan group consists of infinite lower triangular matrices that correspond to certain operators in the space of formal power series. Each such group contains the matrix (generalized Pascal matrix), elements of which are generalized binomial coefficients. Generalized Pascal matrices with non-negative elements form an infinite-dimensional vector space. The paper gives an idea of groups similar to the generalized Riordan groups, but associated with matrices, which in the space of generalized Pascal matrices correspond to the points at infinity; examples of such matrices are the matrix of $q$-binomial coefficients for $q=-1$ and the Pascal triangle modulo $2$. An analog of the Lagrange inversion theorem for these groups is given and the corresponding examples are considered.
\end{abstract}

\section{Introduction}
Matrices that we will consider correspond to operators in the space of formal power series over the field of real numbers. Based on this, we associate the rows and columns of matrices with the generating functions of their elements, i.e., formal power series. Thus, the expression $Aa\left( x \right)=b\left( x \right)$ means that the column vector multiplied by the matrix $A$ has the generating function $a\left( x \right)=\sum\nolimits_{n=0}^{\infty }{{{a}_{n}}}{{x}^{n}}$, resultant column vector has the generating function $b\left( x \right)=\sum\nolimits_{n=0}^{\infty }{{{b}_{n}}{{x}^{n}}}$.  The $n$th coefficient of the series $a\left( x \right)$ denote $\left[ {{x}^{n}} \right]a\left( x \right)$; the $\left( n,m \right)$th element, generating functions of the  $n$th row and $n$th column of the matrix $A$ will be denoted  respectively by ${{A}_{n,m}}={{\left( A \right)}_{n,m}}$,   $\left[ n,\to  \right]A$,    $A{{x}^{n}}$.

The matrix $\left( f\left( x \right),xg\left( x \right) \right)$, $n$th column of which has the generating function \linebreak$f\left( x \right){{x}^{n}}{{g}^{n}}\left( x \right)$, is called a Riordan matrix (Riordan array)  [2,11]. It is a product of two matrices that correspond to the operations of multiplication and composition of series:
$$\left( f\left( x \right),xg\left( x \right) \right)=\left( f\left( x \right),x \right)\left( 1,xg\left( x \right) \right),$$
$$\left( f\left( x \right),x \right)a\left( x \right)=f\left( x \right)a\left( x \right), \qquad\left( 1,xg\left( x \right) \right)a\left( x \right)=a\left( xg\left( x \right) \right),$$
$$\left( f\left( x \right),xg\left( x \right) \right)\left( b\left( x \right),xa\left( x \right) \right)=\left( f\left( x \right)b\left( xg\left( x \right) \right),xg\left( x \right)a\left( xg\left( x \right) \right) \right).$$

The matrices $\left( f\left( x \right),xg\left( x \right) \right)$, ${{f}_{0}}\ne 0$, ${{g}_{0}}\ne 0$, form a group called the Riordan group. Matrices of the form $\left( f\left( x \right),x \right)$ form a normal subgroup called the Appell subgroup, matrices of the form $\left( 1,xg\left( x \right) \right)$ form a subgroup called the Lagrange subgroup. The subgroup of the matrices $\left( g\left( x \right),xg\left( x \right) \right)$ isomorphic to the Lagrange subgroup is called the Bell subgroup.

The matrices
$${{\left| {{e}^{x}} \right|}^{-1}}\left( f\left( x \right),xg\left( x \right) \right)\left| {{e}^{x}} \right|={{\left( f\left( x \right),xg\left( x \right) \right)}_{E}},$$
where $\left| {{e}^{x}} \right|$ is the diagonal matrix: $\left| {{e}^{x}} \right|{{x}^{n}}={{{x}^{n}}}/{n!}\;$, are called exponential Riordan matrices [23,25]. In this connection, the matrices $\left( f\left( x \right),xg\left( x \right) \right)$ are called ordinary Riordan matrices. Denote $\left[ n,\to  \right]{{\left( f\left( x \right),xg\left( x \right) \right)}_{E}}={{s}_{n}}\left( x \right)$. Then
$$\sum\limits_{n=0}^{\infty }{\frac{{{s}_{n}}\left( \varphi  \right)}{n!}{{x}^{n}}}=f\left( x \right)\exp \left( \varphi xg\left( x \right) \right).$$
In the general case (but for ${{f}_{0}},{{g}_{0}}\ne 0$), the sequence of polynomials ${{s}_{n}}\left( x \right)$ is called the Sheffer sequence; in the case $g\left( x \right)=1$ – the Appell sequence, in the case $f\left( x \right)=1$ – the binomial sequence. Properties of Shaeffer sequences are a subject of study of the umbral  calculus [19]

The matrix $P$ whose power  is defined by the identity 
$${{P}^{\varphi }}=\left( \frac{1}{1-\varphi x},\frac{x}{1-\varphi x} \right)={{\left( {{e}^{\varphi x}},x \right)}_{E}}$$ 
is called the Pascal matrix.

The matrices
$${{\left| c\left( x \right) \right|}^{-1}}\left( f\left( x \right),xg\left( x \right) \right)\left| c\left( x \right) \right|={{\left( f\left( x \right),xg\left( x \right) \right)}_{c\left( x \right)}}$$
where $\left| c\left( x \right) \right|$ is the diagonal matrix: $\left| c\left( x \right) \right|{{x}^{n}}={{c}_{n}}{{x}^{n}}$, ${{c}_{n}}\ne 0$, are called generalized Riordan matrices [23] (other options: $\left( c \right)$-Riordan matrices [8], $W$-Riordan matrices [24]). Denote $\left[ n,\to  \right]{{\left( f\left( x \right),xg\left( x \right) \right)}_{c\left( x \right)}}={{\hat{s}}_{n}}\left( x \right)$. Then
$$\sum\limits_{n=0}^{\infty }{{{c}_{n}}{{{\hat{s}}}_{n}}}\left( \varphi  \right){{x}^{n}}=f\left( x \right)c\left( \varphi xg\left( x \right) \right).$$
The sequence ${{\left\{ {{{\hat{s}}}_{n}}\left( x \right) \right\}}_{n\ge 0}}$ is called the generalized Sheffer sequence, the sequence\linebreak ${{\left\{ {{c}_{n}}{{{\hat{s}}}_{n}}\left( x \right) \right\}}_{n\ge 0}}$ is called the Boas-Buck sequence. Properties of these sequences are a subject of study of the non-classical umbral calculus [19]. 

Generalized Riordan matrices are related to the following generalization of binomial coefficients [7]. For the coefficients of the formal power series $b\left( x \right)$, ${{b}_{0}}=0$; ${{b}_{n}}\ne 0$,  $n>0$, denote
$${{b}_{0}}!=1,  \qquad{{b}_{n}}!=\prod\limits_{m=1}^{n}{{{b}_{m}}},   \qquad{{\left( \begin{matrix}
   n  \\
   m  \\
\end{matrix} \right)}_{b}}=\frac{{{b}_{n}}!}{{{b}_{m}}!{{b}_{n-m}}!};   \qquad{{\left( \begin{matrix}
   n  \\
   m  \\
\end{matrix} \right)}_{b}}=0,  \qquad m>n.$$
Then
$${{\left( \begin{matrix}
   n  \\
   m  \\
\end{matrix} \right)}_{b}}={{\left( \begin{matrix}
   n-1  \\
   m-1  \\
\end{matrix} \right)}_{b}}+\frac{{{b}_{n}}-{{b}_{m}}}{{{b}_{n-m}}}{{\left( \begin{matrix}
   n-1  \\
   m  \\
\end{matrix} \right)}_{b}}.$$
Denote
$${{P}_{c\left( x \right)}}={{\left( c\left( x \right),x \right)}_{c\left( x \right)}}, \qquad{{P}_{c\left( x \right)}}x=b\left( x \right).$$
Then
$${{b}_{n}}=\frac{{{c}_{1}}{{c}_{n-1}}}{{{c}_{n}}}, \qquad{{c}_{n}}=\frac{c_{1}^{n}}{{{b}_{n}}!} ,     \qquad{{\left( {{P}_{c\left( x \right)}} \right)}_{n,m}}=\frac{{{c}_{m}}{{c}_{n-m}}}{{{c}_{n}}}={{\left( \begin{matrix}
   n  \\
   m  \\
\end{matrix} \right)}_{b}}.$$

In the general case, for generalized binomial coefficients we will use the notation ${{\left( {{P}_{c\left( x \right)}} \right)}_{n,m}}={n\choose m}_{c\left( x \right)}$. In various specific cases, we will use alternative notation. Since ${{P}_{c\left( x \right)}}={{P}_{c\left( \varphi x \right)}}$, we assume for unambiguity that ${{c}_{1}}=1$. The matrix ${{P}_{c\left( x \right)}}$  will be called the generalized Pascal matrix. To each matrix ${{P}_{c\left( x \right)}}$  we associate the generalized Riordan group, of which it is an element. The series $c\left( x \right)$ will be called the group parameter, the matrix ${{P}_{c\left( x \right)}}$ will be called the central element of the group.

The Pascal matrix provides an example of the ambiguity that arises in the classification of generalized Riordan matrices. It is the ordinary (i.e., no central) element of the ordinary Riordan group and the central element of the exponential Riordan group. The problem caused by this ambiguity is solved in the paper [3]. This paper contains a list of all triplets of the  series
${{f}_{1}}\left( x \right)$, ${{g}_{1}}\left( x \right)$, $c\left( x \right)$, ${{f}_{1}}\left( 0 \right)=1$, ${{g}_{1}}\left( 0 \right)\ne 0$, such that
$${{\left( {{f}_{1}}\left( x \right),x{{g}_{1}}\left( x \right) \right)}_{c\left( x \right)}}=\left( {{f}_{2}}\left( x \right),x{{g}_{2}}\left( x \right) \right) \eqno(1)$$
and, respectively,
$$\left( {{f}_{1}}\left( x \right),x{{g}_{1}}\left( x \right) \right)={{\left( {{f}_{2}}\left( x \right),x{{g}_{2}}\left( x \right) \right)}_{\tilde{c}\left( x \right)}}, \qquad\tilde{c}\left( x \right)=\sum\limits_{n=0}^{\infty }{\frac{{{x}^{n}}}{{{c}_{n}}}}.$$
In [3], equation (1) is presented in a simpler form. Let $D$ be the matrix of the differentiation operator: $D{{x}^{n}}=n{{x}^{n-1}}$. Denote $\left( x,x \right)D=\tilde{D}$, $\tilde{D}{{x}^{n}}=n{{x}^{n}}$. Sinse
$$D\left( f\left( x \right),x \right)=\left( f\left( x \right),x \right)D+\left( {f}'\left( x \right),x \right),\qquad 
D\left( 1,xg\left( x \right) \right)=\left( {{\left( xg\left( x \right) \right)}^{\prime }},xg\left( x \right) \right)D,$$
then $\left( f\left( x \right),xg\left( x \right) \right)\tilde{D}=L\left( f\left( x \right),xg\left( x \right) \right)$, where 
$$L=\left( \frac{g\left( x \right)}{{{\left( xg\left( x \right) \right)}^{\prime }}},x \right)\tilde{D}-\left( \frac{xg\left( x \right){f}'\left( x \right)}{{{\left( xg\left( x \right) \right)}^{\prime }}f\left( x \right)},x \right).$$
Thus, equation (1) is reduced to the equation ${{\left| c\left( x \right) \right|}^{-1}}{{L}_{1}}\left| c\left( x \right) \right|={{L}_{2}}$. 

For the fixed $m\ge 1$  we denote 
$${{c}^{*}}\left( x \right)=\sum\limits_{n=0}^{\infty }{c_{n}^{*}}{{x}^{n}}, \qquad c_{0}^{*}=1, \qquad c_{km+i}^{*}={{\left( c_{m}^{*} \right)}^{k}}c_{i}^{*} , \qquad k\ge 0, \qquad 0\le i<m.$$
All solutions of equation (1) are contained in the next six points, some of which intersect (we will limit ourselves to the condition $g\left( 0 \right)=1$):\\
1.
$$c\left( x \right)={{c}^{*}}\left( x \right),$$
$${{\left( f\left( {{x}^{m}} \right),xg\left( {{x}^{m}} \right) \right)}_{c\left( x \right)}}=\left( f\left( {{{x}^{m}}}/{{{c}_{m}}}\; \right),xg\left( {{{x}^{m}}}/{{{c}_{m}}}\; \right) \right).$$
2. 
$${{c}_{{{n}_{0}}}}=1, \qquad{{c}_{n}}=\varphi c_{n}^{*}, \qquad{{n}_{0}}=0; \qquad{{c}_{{{n}_{0}}}}=\varphi c_{{{n}_{0}}}^{*},  \qquad{{c}_{n}}=c_{n}^{*}, \qquad{{n}_{0}}>0,$$
$${{\left( {{g}^{-{{n}_{0}}}}\left( {{x}^{m}} \right),xg\left( {{x}^{m}} \right) \right)}_{c\left( x \right)}}=\left( {{g}^{-{{n}_{0}}}}\left( {{{x}^{m}}}/{{{c}_{m}}}\; \right),xg\left( {{{x}^{m}}}/{{{c}_{m}}}\; \right) \right),$$
where $g\left( x \right)=1+\sum\nolimits_{n=q}^{\infty }{{{g}_{n}}}{{x}^{n}}$, $q>0$, ${{g}_{q}}\ne 0$, ${{n}_{0}}<qm$. \\
3.
$${{c}_{km+i}}=\left\{ \begin{matrix}
   {{\left( {{c}_{m}} \right)}^{k}}{{c}_{i}}\text{,} & \text{if}\quad i\ne {{i}_{0}}\quad\text{ or}\quad i={{i}_{0}},\text{ }k\le {{k}_{0}},  \\
   {{\left( {{c}_{m}} \right)}^{k-{{k}_{0}}-1}}{{c}_{{{n}_{0}}+m}}, & \text{if}      \quad i={{i}_{0}},\text{ }k>{{k}_{0}},  \\
\end{matrix} \right.$$
where ${{n}_{0}}={{k}_{0}}m+{{i}_{0}}$, ${{k}_{0}}\ge 0$, $0\le {{i}_{0}}<m$, ${{n}_{0}}>0$,
$${{\left( {{\left( 1+\varphi {{x}^{m}} \right)}^{\frac{{{n}_{0}}}{m}}},x{{\left( 1+\varphi {{x}^{m}} \right)}^{-\frac{1}{m}}} \right)}_{c\left( x \right)}}=\left( {{\left( 1+\frac{\varphi }{{{c}_{m}}}{{x}^{m}} \right)}^{\frac{{{n}_{0}}}{m}}},x{{\left( 1+\frac{\varphi }{{{c}_{m}}}{{x}^{m}} \right)}^{-\frac{1}{m}}} \right).$$
4.
$${{c}_{km+i}}=\frac{1}{{{\left( \alpha +{i}/{m}\; \right)}_{k}}}{{\left( \alpha {{c}_{m}} \right)}^{k}}{{c}_{i}},$$
where ${{\left( \alpha  \right)}_{k}}=\alpha \left( \alpha +1 \right)...\left( \alpha +k-1 \right)$, $\alpha \ne -{n}/{m}\;$, $n\ge 0$ ,
$${{\left( \exp \left( \frac{\varphi }{m}{{x}^{m}} \right),x \right)}_{c\left( x \right)}}=\left( {{\left( 1-\frac{\varphi }{m\alpha {{c}_{m}}}{{x}^{m}} \right)}^{-\alpha }},x{{\left( 1-\frac{\varphi }{m\alpha {{c}_{m}}}{{x}^{m}} \right)}^{-\frac{1}{m}}} \right).$$
5.
$${{c}_{km+i}}=\frac{{{\left( \beta +{i}/{m}\; \right)}_{k}}}{{{\left( \alpha +{i}/{m}\; \right)}_{k}}}{{\left( \frac{\alpha {{c}_{m}}}{\beta } \right)}^{k}}{{c}_{i}},$$ 
$\alpha \ne \beta $,  $\alpha ,\beta \ne -{n}/{m}\;$, $n\ge 0$,
$${{\left( {{\left( 1+\varphi {{x}^{m}} \right)}^{-\frac{b}{m}}},x{{\left( 1+\varphi {{x}^{m}} \right)}^{-\frac{1}{m}}} \right)}_{c\left( x \right)}}=\left( {{f}_{2}}\left( x \right),x{{g}_{2}}\left( x \right) \right),$$
$${{f}_{2}}\left( x \right)={{\left( 1+\frac{\varphi }{{{c}_{m}}}\left( 1+\frac{\left( b+1 \right)\left( \beta -\alpha  \right)}{\alpha \left( 1+m\beta  \right)} \right){{x}^{m}} \right)}^{-\frac{b}{m\left( 1+\frac{\left( b+1 \right)\left( \beta -\alpha  \right)}{\alpha \left( 1+m\beta  \right)} \right)}}},$$
$${{g}_{2}}\left( x \right)={{\left( 1+\frac{\varphi }{{{c}_{m}}}\left( 1+\frac{\left( b+1 \right)\left( \beta -\alpha  \right)}{\alpha \left( 1+m\beta  \right)} \right){{x}^{m}} \right)}^{-\frac{1}{m}}},$$
$b\ne -n$, $n\ge 0$. If $b=m\beta $, then
$${{\left( {{\left( 1+\varphi {{x}^{m}} \right)}^{-\beta }},x{{\left( 1+\varphi {{x}^{m}} \right)}^{-\frac{1}{m}}} \right)}_{c\left( x \right)}}=\left( {{\left( 1+\frac{\varphi \beta }{\alpha {{c}_{m}}}{{x}^{m}} \right)}^{-\alpha }},x{{\left( 1+\frac{\varphi \beta }{\alpha {{c}_{m}}}{{x}^{m}} \right)}^{-\frac{1}{m}}} \right).$$
6. 
$${{c}_{km+i}}=\frac{\beta +{i}/{m}\;}{\beta +k+{i}/{m}\;}{{\left( \frac{\beta +1}{\beta }{{c}_{m}} \right)}^{k}}{{c}_{i}},$$
$\beta \ne -{n}/{m}\;$, $n\ge 0$,
$${{\left( {{g}^{m\beta }}\left( {{x}^{m}} \right),xg\left( {{x}^{m}} \right) \right)}_{c\left( x \right)}}=\left( {{\left( xg\left( {{{x}^{m}}}/{{{\beta }^{*}}{{c}_{m}}}\; \right) \right)}^{\prime }}{{g}^{m\beta -1}}\left( {{{x}^{m}}}/{{{\beta }^{*}}{{c}_{m}}}\; \right),xg\left( {{{x}^{m}}}/{{{\beta }^{*}}{{c}_{m}}}\; \right) \right),$$
where ${{\beta }^{*}}={\left( \beta +1 \right)}/{\beta }\;$. If we put ${{c}_{i}}={\beta }/{\left( \beta +{i}/{m}\; \right)}\;$, ${{c}_{m}}={\beta }/{\left( \beta +1 \right)}\;$, then ${{c}_{n}}={m\beta }/{\left( m\beta +n \right)}\;$. Note the intersection with point 5: $g\left( {{x}^{m}} \right)={{\left( 1+\varphi {{x}^{m}} \right)}^{{-1}/{m}\;}}$, $\alpha =\beta +1$. 
 
According to the point 5 of the set of solutions to equation (1),
$${{\left( {{\left( 1-\ x \right)}^{-\beta }},x{{\left( 1-\ x \right)}^{-1}} \right)}_{c\left( x \right)}}=\left( {{\left( 1-\frac{\ \beta }{\alpha }x \right)}^{-\alpha }},x{{\left( 1-\frac{\ \beta }{\alpha }x \right)}^{-1}} \right),$$
$${{c}_{n}}=\frac{{{\left( \beta  \right)}_{n}}{{\alpha }^{n}}}{{{\left( \alpha  \right)}_{n}}{{\beta }^{n}}}, \qquad{{\left( \beta  \right)}_{n}}=\beta \left( \beta +1 \right)...\left( \beta +n-1 \right), \quad\alpha \ne \beta ,  \quad\alpha ,\beta \ne -n.$$
Thus, the Pascal matrix is an element of each generalized Riordan group with the parameter
$c\left( x \right)={{\left( 1+\varphi x \right)}^{{1}/{\varphi }\;}}$, $\varphi =-\left( {1}/{\alpha }\; \right)$ :
$$P={{\left( {{\left( 1+\varphi x \right)}^{\frac{1}{\varphi }}},\frac{x}{1+\varphi x} \right)}_{c\left( x \right)}}={{P}_{c\left( x \right)}}{{\left( 1,\frac{x}{1+\varphi x} \right)}_{c\left( x \right)}}.$$
The exponential Riordan group can be considered as a limit of the set of these groups when the $\varphi $  tending  to zero.

In [4], the set of generalized Pascal matrices is extended by matrices of a special form, which are represented by the following construction. We will consider the set of generalized Pascal matrices as a group with respect to the Hadamard multiplication (denoted by $\times $):
$${{P}_{c\left( x \right)}}\times {{P}_{g\left( x \right)}}={{P}_{c\left( x \right)\times g\left( x \right)}},  \qquad c\left( x \right)\times g\left( x \right)=\sum\limits_{n=0}^{\infty }{{{c}_{n}}{{g}_{n}}}{{x}^{n}}.$$
We introduce the special system of matrices
$${{P}_{\varphi ,q}}={{P}_{q}}\left( \varphi  \right)={{P}_{c\left( \varphi ,q,x \right)}},  \qquad c\left( \varphi ,q,x \right)=\left( \sum\limits_{n=0}^{q-1}{{{x}^{n}}} \right){{\left( 1-\frac{{{x}^{q}}}{\varphi } \right)}^{-1}}, \qquad q>1.$$
Then 
$${{c}_{qn+i}}=\frac{1}{{{\varphi }^{n}}}, \quad 0\le i<q; \qquad {{c}_{qn-i}}=\frac{1}{{{\varphi }^{n-1}}},  \quad 0<i\le q,$$
$$\frac{{{c}_{qm+j}}{{c}_{q\left( n-m \right)+i-j}}}{{{c}_{qn+i}}}=\left\{ \begin{matrix}
   \frac{{{\varphi }^{n}}}{{{\varphi }^{m}}{{\varphi }^{n-m}}}=1,\text{  }i\ge j,  \\
   \frac{{{\varphi }^{n}}}{{{\varphi }^{m}}{{\varphi }^{n-m-1}}}=\varphi ,\text{ }i<j,  \\
\end{matrix} \right.$$
or 
$${{\left( {{P}_{\varphi ,q}} \right)}_{n,m}}=\left\{ \begin{matrix}
   1,n\left( \bmod q \right)\ge m\left( \bmod q \right),  \\
   \varphi ,n\left( \bmod q \right)<m\left( \bmod q \right).  \\
\end{matrix} \right.$$
For example, ${{P}_{\varphi ,2}}$, ${{P}_{\varphi ,3}}$:
$$\left( \begin{matrix}
   1 & 0 & 0 & 0 & 0 & 0 & 0 & 0 & 0 & \ldots   \\
   1 & 1 & 0 & 0 & 0 & 0 & 0 & 0 & 0 & \ldots   \\
   1 & \varphi  & 1 & 0 & 0 & 0 & 0 & 0 & 0 & \ldots   \\
   1 & 1 & 1 & 1 & 0 & 0 & 0 & 0 & 0 & \ldots   \\
   1 & \varphi  & 1 & \varphi  & 1 & 0 & 0 & 0 & 0 & \ldots   \\
   1 & 1 & 1 & 1 & 1 & 1 & 0 & 0 & 0 & \dots   \\
   1 & \varphi  & 1 & \varphi  & 1 & \varphi  & 1 & 0 & 0 & \ldots   \\
   1 & 1 & 1 & 1 & 1 & 1 & 1 & 1 & 0 & \ldots   \\
   1 & \varphi  & 1 & \varphi  & 1 & \varphi  & 1 & \varphi  & 1 & \ldots   \\
   \vdots  & \vdots  & \vdots  & \vdots  & \vdots  & \vdots  & \vdots  & \vdots  & \vdots  & \ddots   \\
\end{matrix} \right), \quad \left( \begin{matrix}
   1 & 0 & 0 & 0 & 0 & 0 & 0 & 0 & 0 & \ldots   \\
   1 & 1 & 0 & 0 & 0 & 0 & 0 & 0 & 0 & \ldots   \\
   1 & 1 & 1 & 0 & 0 & 0 & 0 & 0 & 0 & \ldots   \\
   1 & \varphi  & \varphi  & 1 & 0 & 0 & 0 & 0 & 0 & \ldots   \\
   1 & 1 & \varphi  & 1 & 1 & 0 & 0 & 0 & 0 & \ldots   \\
   1 & 1 & 1 & 1 & 1 & 1 & 0 & 0 & 0 & \ldots   \\
   1 & \varphi  & \varphi  & 1 & \varphi  & \varphi  & 1 & 0 & 0 & \ldots   \\
   1 & 1 & \varphi  & 1 & 1 & \varphi  & 1 & 1 & 0 & \ldots   \\
   1 & 1 & 1 & 1 & 1 & 1 & 1 & 1 & 1 & \ldots   \\
   \vdots  & \vdots  & \vdots  & \vdots  & \vdots  & \vdots  & \vdots  & \vdots  & \vdots  & \ddots   \\
\end{matrix} \right).$$
The matrices $\left( a\left( x \right)|{{P}_{\varphi ,q}}\times {{P}_{c\left( x \right)}} \right)$, where $\varphi $, $q$ and $c\left( x \right)$ are fixed, $a\left( x \right)\in \mathbb{R}\left[ \left[ x \right] \right]$,\linebreak  ${{\left( a\left( x \right)|{{P}_{\varphi ,q}}\times {{P}_{c\left( x \right)}} \right)}_{n,m}}={{a}_{n-m}}{{\left( {{P}_{\varphi ,q}}\times {{P}_{c\left( x \right)}} \right)}_{n,m}}$, form an algebra isomorphic to the algebra of formal power series with the multiplication
$$a\left( x \right)\circ b\left( x \right)=\sum\limits_{n=0}^{\infty }{\left( \sum\limits_{m=0}^{n}{{{\left( {{P}_{\varphi ,q}}\times {{P}_{c\left( x \right)}} \right)}_{n,m}}{{a}_{n-m}}}{{b}_{m}} \right){{x}^{n}}}:$$
$$\left( a\left( x \right)|{{P}_{\varphi ,q}}\times {{P}_{c\left( x \right)}} \right)\left( b\left( x \right)|{{P}_{\varphi ,q}}\times {{P}_{c\left( x \right)}} \right)=\left( a\left( x \right)\circ b\left( x \right)|{{P}_{\varphi ,q}}\times {{P}_{c\left( x \right)}} \right).$$
This proposition remains true for any values of $\varphi$, i.e., the case $\varphi = 0$ is not exceptional. The matrix ${{P}_{0,q}}\times {{P}_{c\left( x \right)}}$, as well as the Hadamard product of such matrices, will be called the zero generalized Pascal matrix.

 The zero generalized Pascal matrix appears when considering the set of matrices ${{P}_{g\left( q,x \right)}}$:
$${{P}_{g\left( q,x \right)}}{{x}^{n}}={{x}^{n}}\prod\limits_{m=0}^{n}{{{\left( 1-{{q}^{m}}x \right)}^{-1}}}, \qquad q\in \mathbb{R}.$$
Here $g\left( 0,x \right)={{\left( 1-x \right)}^{-1}}$,  $g\left( 1,x \right)={{e}^{x}}$. In other cases, except $q=-1$,
$$g\left( q,x \right)=\sum\limits_{n=0}^{\infty }{\frac{{{\left( q-1 \right)}^{n}}}{\left( {{q}^{n}}-1 \right)!}}{{x}^{n}},  \qquad\left( {{q}^{n}}-1 \right)!=\prod\limits_{m=1}^{n}{\left( {{q}^{m}}-1 \right)},  \qquad\left( {{q}^{0}}-1 \right)!=1.$$
For the matrix 
$${{P}_{g\left( -1,x \right)}}=\left( \begin{matrix}
   1 & 0 & 0 & 0 & 0 & 0 & 0 & \ldots   \\
   1 & 1 & 0 & 0 & 0 & 0 & 0 & \ldots   \\
   1 & 0 & 1 & 0 & 0 & 0 & 0 & \ldots   \\
   1 & 1 & 1 & 1 & 0 & 0 & 0 & \ldots   \\
   1 & 0 & 2 & 0 & 1 & 0 & 0 & \ldots   \\
   1 & 1 & 2 & 2 & 1 & 1 & 0 & \ldots   \\
   1 & 0 & 3 & 0 & 3 & 0 & 1 & \ldots   \\
   \vdots  & \vdots  & \vdots  & \vdots  & \vdots  & \vdots  & \vdots  & \ddots   \\
\end{matrix} \right)$$
the series $g\left( -1,x \right)$ is not defined, but since
$${{\left( {{P}_{g\left( -1,x \right)}} \right)}_{2n+i,2m+j}}=\left[ {{x}^{2n+i}} \right]\frac{{{\left( 1+x \right)}^{1-j}}{{x}^{2m+j}}}{{{\left( 1-{{x}^{2}} \right)}^{m+1}}}=\left\{ \begin{matrix}
   \left( \begin{matrix}
   n  \\
   m  \\
\end{matrix} \right),i\ge j,  \\
   0,i<j,  \\
\end{matrix} \right. \qquad i,j=0,1,$$
 then
$${{P}_{g\left( -1,x \right)}}=P{}_{0,2}\times {{P}_{c\left( x \right)}},  \qquad c\left( x \right)=\left( 1+x \right){{e}^{{{x}^{2}}}}:$$
$${{c}_{2n+i}}=\frac{1}{n!},  \qquad 0\le i<2;  \qquad{{c}_{2n-i}}=\frac{1}{\left( n-1 \right)!},\qquad 0<i\le 2,$$
$${{\left( {{P}_{c\left( x \right)}} \right)}_{2n+i,2m+j}}=\left\{ \begin{matrix}
   \left( \begin{matrix}
   n  \\
   m  \\
\end{matrix} \right),i\ge j,  \\
   n\left( \begin{matrix}
   n-1  \\
   m  \\
\end{matrix} \right),i<j.  \\
\end{matrix} \right.$$

Note that the matrix  ${{P}_{g\left( -1,x \right)}}$ is the A051159, [20]. In [10], this matrix is called the Pauli Pascal triangle for the following reason. Let ${{\sigma }_{x}}$, ${{\sigma }_{y}}$ be the Pauli spin matrices:
 $\sigma _{x}^{2}=\sigma _{y}^{2}=1$, ${{\sigma }_{x}}{{\sigma }_{y}}=-{{\sigma }_{y}}{{\sigma }_{x}}$. Then
$${{\left( {{\sigma }_{x}}+{{\sigma }_{y}} \right)}^{n}}=\sum\limits_{m=0}^{n}{{{\left( {{P}_{g\left( -1,x \right)}} \right)}_{n,m}}}\sigma _{x}^{n-m}\sigma _{y}^{m}.$$

Each "nonzero" generalized Pascal matrix can be expanded in the Hadamard product of the matrices ${{P}_{\varphi ,q}}$ (the expansion reduces to the  expansion  of the first column of the matrix ${{P}_{c\left( x \right)}}$  in the Hadamard product of the first columns of the matrices ${{P}_{\varphi ,q}}$). Let  ${{e}_{q}}$ be the basis vector of the infinite-dimensional vector space over the field of real numbers, the numbering of the basis vectors of which begins with two. The mapping  ${{P}_{\varphi ,q}}\to {{e}_{q}}\log \left| \varphi  \right|$ of the set of generalized Pascal matrices to the infinite-dimensional vector space is a group homomorphism whose kernel  consists of all involutions of the group of generalized Pascal matrices, i.e., of matrices whose nonzero elements are $\pm 1$. Thus, the set of generalized Pascal matrices whose elements are non-negative numbers is an infinite-dimensional vector space. Zero generalized Pascal matrices can be considered as points at infinity in this space. 

The second special system of generalized Pascal matrices ${{P}_{\left[ \varphi ,q \right]}}$ is related to the system of matrices  ${{P}_{\varphi ,q}}$  as follows:
$${{P}_{\left[ \varphi ,q \right]}}={{P}_{\varphi ,q}}\times {{P}_{\varphi ,{{q}^{2}}}}\times {{P}_{\varphi ,{{q}^{3}}}}\times ...=\prod\limits_{k=1}^{\infty }{\times {{P}_{\varphi ,{{q}^{k}}}}}.$$
In [4], the matrices ${{P}_{\left[ \varphi ,q \right]}}$ are called fractal generalized Pascal matrices because of their fractal properties. Denote 
 ${{P}_{\left[ q,q \right]}}={{P}_{\left[ q \right]}}$. The series ${{P}_{\left[ q \right]}}x$ is a generating function of the distribution of the divisors ${{q}^{k}}$ in the series of natural numbers. For example,
$${{P}_{\left[ 2 \right]}}x=x+2{{x}^{2}}+{{x}^{3}}+4{{x}^{4}}+{{x}^{5}}+2{{x}^{6}}+{{x}^{7}}+8{{x}^{8}}+{{x}^{9}}+2{{x}^{10}}+{{x}^{11}}+4{{x}^{12}}+...,$$
$${{P}_{\left[ 3 \right]}}x=x+{{x}^{2}}+3{{x}^{3}}+{{x}^{4}}+{{x}^{5}}+3{{x}^{6}}+{{x}^{7}}+{{x}^{8}}+9{{x}^{9}}+{{x}^{10}}+{{x}^{11}}+3{{x}^{12}}+....$$
Since the series $Px$  is a generating function of the sequence of natural numbers, then
$$P=\prod\limits_{k=1}^{\infty }{\times {{P}_{\left[ {{p}_{k}} \right]}}},$$
where ${{\left\{ {{p}_{k}} \right\}}_{k\ge 1}}$ is the sequence of prime numbers. Note that ${{P}_{\left[ q \right]}}={{P}_{c\left( x \right)}}$,
$$c\left( x \right)=\left( \sum\limits_{n=0}^{q-1}{{{x}^{n}}} \right)c\left( \frac{{{x}^{q}}}{q} \right)=\prod\limits_{n=0}^{\infty }{{{w}_{q}}}\left( {{{x}^{{{q}^{n}}}}}/{{{q}^{\frac{{{q}^{n}}-1}{q-1}}}}\; \right), \qquad{{w}_{q}}\left( x \right)=\left( \sum\limits_{n=0}^{q-1}{{{x}^{n}}} \right).$$

The matrices ${{P}_{\left[ b \right]}}$, under a general name the $b$-binomial triangle, were introduced in [1] in connection with the generalization of theorems on a divisibility of binomial coefficients.

The purpose of this  paper is to extend the set of generalized Riordan groups at the expense of groups associated with zero generalized Pascal matrices. A zero generalized Pascal matrix of general form will be denoted by $ {{P} _ {0}} $. Section 2 gives an idea of the matrix algebra whose element is the matrix $ {{P} _ {0}} $. We denote this algebra by $\left[\!\left[ {{P}_{0}} \right]\!\right]$. In sections 3 and 4, as examples of constructions characteristic for the algebra $\left[\!\left[ {{P}_{0}} \right]\!\right]$, we consider the "block groups" (algebra $\left[\!\left[ {{P}_{0,q}} \right]\!\right]$) and the "fractal groups" (algebra $\left[\!\left[ {{P}_{\left[ 0,q \right]}} \right]\!\right]$). Section 5 gives an idea of the group $R\left( {{P}_{0}} \right)$ similar to the generalized Riordan group $R\left( {{P}_{c\left( x \right)}} \right)$. An analog of the Lagrange inversion theorem for this group is given and corresponding examples are considered.
\section{Algebra $\left[\!\left[ {{P}_{0}} \right]\!\right]$}
 We denote
$${{\left( {{P}_{0}} \right)}_{n,m}}={{\left( \begin{matrix}
   n  \\
   m  \\
\end{matrix} \right)}_{0}}, \qquad{{\left( a\left( x \right)|{{P}_{0}} \right)}_{n,m}}={{a}_{n-m}}{{\left( \begin{matrix}
   n  \\
   m  \\
\end{matrix} \right)}_{0}},$$
$$\left( a\left( x \right)|{{P}_{0}} \right)\left( b\left( x \right)|{{P}_{0}} \right)=\left( a\left( x \right)\circ b\left( x \right)|{{P}_{0}} \right), \qquad{{\left( a\left( x \right)|{{P}_{0}} \right)}^{n}}=\left( {{a}^{\left( n \right)}}\left( x \right)|{{P}_{0}} \right).$$
For the case ${{a}_{0}}=1$, we define the power and the logarithm of the matrix 
 $\left( a\left( x \right)|{{P}_{0}} \right)$ in the standard way:
$${{\left( a\left( x \right)|{{P}_{0}} \right)}^{\varphi }}=\sum\limits_{n=0}^{\infty }{\left( \begin{matrix}
   \varphi   \\
   n  \\
\end{matrix} \right)}{{\left( a\left( x \right)-1|{{P}_{0}} \right)}^{n}}=\left( {{a}^{\left( \varphi  \right)}}\left( x \right)|{{P}_{0}} \right),$$
$$\log \left( a\left( x \right)|{{P}_{0}} \right)=\sum\limits_{n=1}^{\infty }{\frac{{{\left( -1 \right)}^{n-1}}}{n}}{{\left( a\left( x \right)-1|{{P}_{0}} \right)}^{n}}=\left( \log \circ a\left( x \right)|{{P}_{0}} \right).$$
For the case ${{b}_{0}}=0$, we define the exponent:
$$\exp \left( b\left( x \right)|{{P}_{0}} \right)=\sum\limits_{n=0}^{\infty }{\frac{{{\left( b\left( x \right)|{{P}_{0}} \right)}^{n}}}{n!}}=\left( \exp \circ b\left( x \right)|{{P}_{0}} \right).$$
Then
$$\left( {{a}^{\left( \varphi  \right)}}\left( x \right)|{{P}_{0}} \right)=\sum\limits_{n=0}^{\infty }{\frac{{{\varphi }^{n}}}{n!}}{{\left( \log \circ a\left( x \right)|{{P}_{0}} \right)}^{n}}.$$

The algebra of matrices $\left( a\left( x \right)|{{P}_{0}} \right)$ will be denoted by $\left[\!\left[ {{P}_{0}} \right]\!\right]$.
The algebra of formal power series isomorphic to the algebra $\left[\!\left[ {{P}_{0}} \right]\!\right]$ will be denoted by 
$\left[\!\left[ {{P}_{0}},a\left( x \right) \right]\!\right]$. For this algebra, we have:
$${{a}^{\left( \varphi  \right)}}\left( x \right)=\sum\limits_{n=0}^{\infty }{\left( \begin{matrix}
   \varphi   \\
   n  \\
\end{matrix} \right){{\left( a\left( x \right)-1 \right)}^{\left( n \right)}}},  \qquad\log \circ a\left( x \right)=\sum\limits_{n=1}^{\infty }{\frac{{{\left( -1 \right)}^{n-1}}}{n}}{{\left( a\left( x \right)-1 \right)}^{\left( n \right)}},$$
$$\exp \circ b\left( x \right)=\sum\limits_{n=0}^{\infty }{\frac{{{b}^{\left( n \right)}}\left( x \right)}{n!}},  \qquad{{a}^{\left( \varphi  \right)}}\left( x \right)=\sum\limits_{n=0}^{\infty }{\frac{{{\varphi }^{n}}}{n!}}{{\left( \log \circ a\left( x \right) \right)}^{\left( n \right)}}=\sum\limits_{n=0}^{\infty }{{{c}_{n}}\left( \varphi  \right){{x}^{n}}},$$
where ${{c}_{n}}\left( \varphi  \right)$ are polynomials in $\varphi $ of degree $\le n$ similar to convolution polynomials [13],
$${{a}^{\left( \varphi  \right)}}\left( x \right)\circ {{a}^{\left( \beta  \right)}}\left( x \right)={{a}^{\left( \varphi +\beta  \right)}}\left( x \right),  \qquad{{a}^{\left( \varphi  \right)}}\left( x \right)\circ {{b}^{\left( \varphi  \right)}}\left( x \right)={{\left( a\left( x \right)\circ b\left( x \right) \right)}^{\left( \varphi  \right)}},$$
$$\log \circ {{a}^{\left( \varphi  \right)}}\left( x \right)=\varphi \log \circ a\left( x \right), \qquad\log \circ \left( a\left( x \right)\circ b\left( x \right) \right)=\log \circ a\left( x \right)+\log \circ b\left( x \right).$$

Let $D$ be the matrix of the differentiation operator. Since $\left( x,x \right)D$  is a diagonal matrix, the identity  
$$\left( x,x \right)D\left( a\left( x \right),x \right)=\left( a\left( x \right),x \right)\left( x,x \right)D+\left( x{a}'\left( x \right),x \right)$$
implies the identity
	$$\left( x,x \right)D\left( a\left( x \right)|{{P}_{0}} \right)=\left( a\left( x \right)|{{P}_{0}} \right)\left( x,x \right)D+\left( x{a}'\left( x \right)|{{P}_{0}} \right).$$ 
For the algebra $\left[\!\left[ {{P}_{0}},a\left( x \right) \right]\!\right]$, we have: 
$$x{{\left( a\left( x \right)\circ b\left( x \right) \right)}^{\prime }}=a\left( x \right)\circ x{b}'\left( x \right)+x{a}'\left( x \right)\circ b\left( x \right),$$
$$x{{\left( {{a}^{\left( n \right)}}\left( x \right) \right)}^{\prime }}=n{{a}^{\left( n-1 \right)}}\left( x \right)\circ x{a}'\left( x \right),$$ 
$$x{{\left( {{a}^{\left( \varphi  \right)}}\left( x \right) \right)}^{\prime }}=x{a}'\left( x \right)\circ \varphi \sum\limits_{n=1}^{\infty }{\left( \begin{matrix}
   \varphi -1  \\
   n-1  \\
\end{matrix} \right){{\left( a\left( x \right)-1 \right)}^{\left( n-1 \right)}}}=\varphi {{a}^{\left( \varphi -1 \right)}}\left( x \right)\circ x{a}'\left( x \right),$$
$$x{{\left( \log \circ a\left( x \right) \right)}^{\prime }}=x{a}'\left( x \right)\circ \sum\limits_{n=1}^{\infty }{{{\left( -1 \right)}^{n-1}}{{\left( a\left( x \right)-1 \right)}^{\left( n-1 \right)}}}=x{a}'\left( x \right)\circ {{a}^{\left( -1 \right)}}\left( x \right).$$

Let ${{P}_{0}}={{P}_{0,q}}$, ${n\choose m}_{0}={n\choose m}_{0,q}$.\\
{\bfseries Theorem 2.1.} \emph{The matrices $\left( {{x}^{qp+i}}|P{}_{0,q} \right)$, $q-\left\lfloor {q}/{2}\; \right\rfloor \le i<q$, form a closed system of zero divisors, i.e., their products with each other and with themselves are equal to zero. }\\ 
{\bfseries Proof.} Since
$$\left[ {{x}^{n}} \right]a\left( x \right)\circ b\left( x \right)=\sum\limits_{m=0}^{n}{{{\left( \begin{matrix}
   n  \\
   m  \\
\end{matrix} \right)}_{0,q}}{{a}_{m}}{{b}_{n-m}}},$$
then
$${{x}^{q{{p}_{1}}+i}}\circ {{x}^{q{{p}_{2}}+j}}={{\left( \begin{matrix}
   q\left( {{p}_{1}}+{{p}_{2}} \right)+i+j  \\
   q{{p}_{1}}+i  \\
\end{matrix} \right)}_{0,q}}{{x}^{q\left( {{p}_{1}}+{{p}_{2}} \right)+i+j}},\quad 0\le i,j<q.$$
Let $i=q-\tilde{i}$, $j=q-\tilde{j}$, $0<\tilde{i},\tilde{j}\le \left\lfloor {q}/{2}\; \right\rfloor $. Then $\left( i+j \right)\left( \bmod q \right)=q-\tilde{i}-\tilde{j}$. Thus, if $q-\left\lfloor {q}/{2}\; \right\rfloor \le i,j<q$, then $\left( i+j \right)\left( \bmod q \right)<i$, or
$$\left( q\left( {{p}_{1}}+{{p}_{2}} \right)+i+j \right)\left( \bmod q \right)<\left( q{{p}_{1}}+i \right)\left( \bmod q \right),$$
$${{\left( \begin{matrix}
   q\left( {{p}_{1}}+{{p}_{2}} \right)+i+j  \\
   q{{p}_{1}}+i  \\
\end{matrix} \right)}_{0,q}}=0.  \qquad\square $$ 

Thus, the series of the form
$${{\eta }_{q}}\left( x \right)=\sum\limits_{n=0}^{\infty }{\sum\limits_{i=q-\left\lfloor {q}/{2}\; \right\rfloor }^{q-1}{{{\eta }_{qn+i}}}}{{x}^{qn+i}}\eqno (2)$$ 
is a nilpotent of degree 2. The series of the form ${{\omega }_{q}}\left( x \right)=1+{{\eta }_{q}}\left( x \right)$  is an unipotent and can be represented as ${{\omega }_{q}}\left( x \right)=1+\log \circ {{\omega }_{q}}\left( x \right)$. These series form a group whose elements are multiplied by the rule
$${{\omega }_{q,1}}\left( x \right)\circ {{\omega }_{q,2}}\left( x \right)=1+\log \circ \left( {{\omega }_{q,1}}\left( x \right)\circ {{\omega }_{q,2}}\left( x \right) \right).$$
\section{Block groups}
Block matrices are an attribute of the algebras associated with the matrix ${{P}_{0,q}}$. Denote ${{\left[ b\left( x \right) \right]}_{q}}=\sum\nolimits_{n=0}^{q-1}{{{b}_{n}}}{{x}^{n}}$. The matrix $\left( {{\left[ b\left( x \right) \right]}_{q}}a\left( {{x}^{q}} \right)|{{P}_{0,q}} \right)$ is a block matrix, the $\left( n,m \right)$th block of which is the matrix ${{a}
_{n-m}}{{\left( b\left( x \right),x \right)}_{q}}$, where ${{\left( b\left( x \right),x \right)}_{q}}$ is the matrix consisting of the first $q$ rows of the matrix $\left( b\left( x \right),x \right)$. For example,
$$\left( {{\left[ b\left( x \right) \right]}_{3}}a\left( {{x}^{3}} \right)|{{P}_{0,3}} \right)=\left( \begin{matrix}
   {{a}_{0}}{{b}_{0}} & 0 & 0 & 0 & 0 & 0 & \ldots   \\
   {{a}_{0}}{{b}_{1}} & {{a}_{0}}{{b}_{0}} & 0 & 0 & 0 & 0 & \ldots   \\
   {{a}_{0}}{{b}_{2}} & {{a}_{0}}{{b}_{1}} & {{a}_{0}}{{b}_{0}} & 0 & 0 & 0 & \ldots   \\
   {{a}_{1}}{{b}_{0}} & 0 & 0 & {{a}_{0}}{{b}_{0}} & 0 & 0 & \ldots   \\
   {{a}_{1}}{{b}_{1}} & {{a}_{1}}{{b}_{0}} & 0 & {{a}_{0}}{{b}_{1}} & {{a}_{0}}{{b}_{0}} & 0 & \ldots   \\
   {{a}_{1}}{{b}_{2}} & {{a}_{1}}{{b}_{1}} & {{a}_{1}}{{b}_{0}} & {{a}_{0}}{{b}_{2}} & {{a}_{0}}{{b}_{1}} & {{a}_{0}}{{b}_{0}} & \ldots   \\
   \vdots  & \vdots  & \vdots  & \vdots  & \vdots  & \vdots  & \ddots   \\
\end{matrix} \right).$$
Hence,  if ${{a}_{0}}\ne 0$, ${{b}_{0}}\ne 0$, these matrices form a group whose elements are multiplied by the rule
$$\left( {{\left[ {{b}_{1}}\left( x \right) \right]}_{q}}{{a}_{1}}\left( {{x}^{q}} \right)|{{P}_{0,q}} \right)\left( {{\left[ {{b}_{2}}\left( x \right) \right]}_{q}}{{a}_{2}}\left( {{x}^{q}} \right)|{{P}_{0,q}} \right)=
\left( {{\left[ {{b}_{1}}\left( x \right){{b}_{2}}\left( x \right) \right]}_{q}}{{a}_{1}}\left( {{x}^{q}} \right){{a}_{2}}\left( {{x}^{q}} \right)|{{P}_{0,q}} \right).$$
Note the intersection of this group with the Riordan group:
$\left( a\left( {{x}^{q}} \right)|{{P}_{0,q}} \right)=\left( a\left( {{x}^{q}} \right),x \right)$.
{\bfseries Theorem 3.1.} 
$$\log \left( {{\left[ b\left( x \right) \right]}_{q}}a\left( {{x}^{q}} \right)|{{P}_{0,q}} \right)=\left( {{\left[ \log b\left( x \right) \right]}_{q}}+\log a\left( {{x}^{q}} \right)|{{P}_{0,q}} \right).$$
{\bfseries Proof.}  It follows from the identities
$$\left( {{\left[ b\left( x \right) \right]}_{q}}a\left( {{x}^{q}} \right)|{{P}_{0,q}} \right)=\left( {{\left[ b\left( x \right) \right]}_{q}}|{{P}_{0,q}} \right)\left( a\left( {{x}^{q}} \right)|{{P}_{0,q}} \right),$$
$$\log {{\left( b\left( x \right),x \right)}_{q}}={{\left( \log b\left( x \right),x \right)}_{q}},\qquad  \log \left( a\left( {{x}^{q}} \right),x \right)=\left( \log a\left( {{x}^{q}} \right),x \right). \qquad\square $$

Thus, since
$${{P}_{0,q}}=\left( \frac{1}{1-x}|{{P}_{0,q}} \right)=\left( {{\left[ \frac{1}{1-x} \right]}_{q}}\frac{1}{1-{{x}^{q}}}|{{P}_{0,q}} \right),$$
then
$$\log {{P}_{0,q}}=\left( \sum\limits_{m=1}^{q-1}{\frac{{{x}^{m}}}{m}+\sum\limits_{m=1}^{\infty }{\frac{{{x}^{mq}}}{m}}|{{P}_{0,q}}} \right).$$
Note that
$${{\left( {{\left[ b\left( x \right) \right]}_{q}}a\left( {{x}^{q}} \right)|{{P}_{0,q}} \right)}^{\varphi }}=\left( {{\left[ {{b}^{\varphi }}\left( x \right) \right]}_{q}}{{a}^{\varphi }}\left( {{x}^{q}} \right)|{{P}_{0,q}} \right).\eqno (3)$$

A generalization of the matrix ${{P}_{g\left( -1,x \right)}}$ considered in the introduction is the matrix  ${{P}_{0,{{c}_{q}}\left( x \right)}}$,  associated with the matrix ${{P}_{c\left( x \right)}}$ as follows:
$${{P}_{0,{{c}_{q}}\left( x \right)}}={{P}_{0,q}}\times {{P}_{{{c}_{q}}\left( x \right)}},  \qquad{{c}_{q}}\left( x \right)={{\left[ {{\left( 1-x \right)}^{-1}} \right]}_{q}}c\left( {{x}^{q}} \right),$$  
$$\left[ {{x}^{qn+i}} \right]{{c}_{q}}\left( x \right)={{c}_{n}},    \qquad 0\le i<q;  \qquad\left[ {{x}^{qn-i}} \right]{{c}_{q}}\left( x \right)={{c}_{n-1}},  \qquad 0<i\le q,$$
$${{\left( {{P}_{{{c}_{q}}\left( x \right)}} \right)}_{qn+i,qm+j}}=\left\{ \begin{matrix}
   {{\left( \begin{matrix}
   n  \\
   m  \\
\end{matrix} \right)}_{c\left( x \right)}},i\ge j,  \\
   {{\left( \begin{matrix}
   n-1  \\
   m  \\
\end{matrix} \right)}_{c\left( x \right)}}{{\left( \begin{matrix}
   n  \\
   1  \\
\end{matrix} \right)}_{c\left( x \right)}},i<j.  \\
\end{matrix} \right.$$

For the algebra of the matrices $\left( a\left( x \right)|{{P}_{c\left( x \right)}} \right)$, ${{\left( a\left( x \right)|{{P}_{c\left( x \right)}} \right)}_{n,m}}={{a}_{n-m}}{n\choose m}_{c\left( x \right)}$, we denote
$$\left( a\left( x \right)|{{P}_{c\left( x \right)}} \right)\left( b\left( x \right)|{{P}_{c\left( x \right)}} \right)=\left( a\left( x \right)*b\left( x \right)|{{P}_{c\left( x \right)}} \right),\qquad\log \left( a\left( x \right)|{{P}_{c\left( x \right)}} \right)=\left( \log *a\left( x \right)|{{P}_{c\left( x \right)}} \right).$$
 Since $\left( a\left( x \right)|{{P}_{c\left( x \right)}} \right)={{\left( a\left( c,x \right),x \right)}_{c\left( x \right)}}$, $a\left( c,x \right)=\left| c\left( x \right) \right|a\left( x \right)$, then 
$$a\left( x \right)*b\left( x \right)={{\left| c\left( x \right) \right|}^{-1}}a\left( c,x \right)b\left( c,x \right),  \qquad\log *a\left( x \right)={{\left| c\left( x \right) \right|}^{-1}}\log a\left( c,x \right).$$
  
The matrix $\left( {{\left[ b\left( x \right) \right]}_{q}}a\left( {{x}^{q}} \right)|{{P}_{0,{{c}_{q}}\left( x \right)}} \right)$ is a block matrix, the $\left( n,m \right)$th block of which is the matrix
$${{a}_{n-m}}{{\left( \begin{matrix}
   n  \\
   m  \\
\end{matrix} \right)}_{c\left( x \right)}}{{\left( b\left( x \right),x \right)}_{q}}.$$
Hence,  if ${{a}_{0}}\ne 0$, ${{b}_{0}}\ne 0$, these matrices form a group whose elements are multiplied by the rule
$$\left( {{\left[ {{b}_{1}}\left( x \right) \right]}_{q}}{{a}_{1}}\left( {{x}^{q}} \right)|{{P}_{0,{{c}_{q}}\left( x \right)}} \right)\left( {{\left[ {{b}_{2}}\left( x \right) \right]}_{q}}{{a}_{2}}\left( {{x}^{q}} \right)|{{P}_{0,{{c}_{q}}\left( x \right)}} \right)=$$
$$=\left( {{\left[ {{b}_{1}}\left( x \right){{b}_{2}}\left( x \right) \right]}_{q}}\left( {{a}_{1}}\left( {{x}^{q}} \right)*{{a}_{2}}\left( {{x}^{q}} \right) \right)|{{P}_{0,{{c}_{q}}\left( x \right)}} \right).$$
{\bfseries Theorem 3.2.}\emph{$$\log \left( {{\left[ b\left( x \right) \right]}_{q}}a\left( {{x}^{q}} \right)|{{P}_{0,{{c}_{q}}\left( x \right)}} \right)=\left( {{\left[ \log b\left( x \right) \right]}_{q}}+{{\left| {{c}_{q}}\left( x \right) \right|}^{-1}}\log a\left( c,{{x}^{q}} \right)|{{P}_{0,{{c}_{q}}\left( x \right)}} \right),$$
where $a\left( c,{{x}^{q}} \right)=\sum\nolimits_{n=0}^{\infty }{{{a}_{n}}{{c}_{n}}{{x}^{qn}}}$}\\
{\bfseries Proof.} Since
$$\left( {{\left[ b\left( x \right) \right]}_{q}}a\left( {{x}^{q}} \right)|{{P}_{0,{{c}_{q}}\left( x \right)}} \right)=\left( {{\left[ b\left( x \right) \right]}_{q}}|{{P}_{0,{{c}_{q}}\left( x \right)}} \right)\left( a\left( {{x}^{q}} \right)|{{P}_{0,{{c}_{q}}\left( x \right)}} \right),$$
$$\left( {{\left[ b\left( x \right) \right]}_{q}}|{{P}_{0,{{c}_{q}}\left( x \right)}} \right)=\left( {{\left[ b\left( x \right) \right]}_{q}}|{{P}_{0,q}} \right),  \quad\left( a\left( {{x}^{q}} \right)|{{P}_{0,{{c}_{q}}\left( x \right)}} \right)={{\left( a\left( c,{{x}^{q}} \right),x \right)}_{{{c}_{q}}\left( x \right)}},$$ 
then
$$\log \circ \left( {{\left[ b\left( x \right) \right]}_{q}}a\left( {{x}^{q}} \right) \right)={{\left[ \log b\left( x \right) \right]}_{q}}+{{\left| {{c}_{q}}\left( x \right) \right|}^{-1}}\log a\left( c,{{x}^{q}} \right). \qquad\square $$ 

Thus, since
$${{P}_{0,{{c}_{q}}\left( x \right)}}=\left( {{\left[ \frac{1}{1-x} \right]}_{q}}\frac{1}{1-{{x}^{q}}}|{{P}_{0,{{c}_{q}}\left( x \right)}} \right),$$
then
$$\log {{P}_{0,{{c}_{q}}\left( x \right)}}=\left( \sum\limits_{m=1}^{q-1}{\frac{{{x}^{m}}}{m}+{{\left| {{c}_{q}}\left( x \right) \right|}^{-1}}\log c\left( {{x}^{q}} \right)|{{P}_{0,{{c}_{q}}\left( x \right)}}} \right).$$
In particular, if $c\left( x \right)={{e}^{x}}$, as in the case of the matrix ${{P}_{g\left( -1,x \right)}}={{P}_{0,{{c}_{2}}\left( x \right)}}$, then
$$\log {{P}_{0,{{c}_{q}}\left( x \right)}}=\left( \sum\limits_{m=1}^{q-1}{\frac{{{x}^{m}}}{m}+{{x}^{q}}|{{P}_{0,{{c}_{q}}\left( x \right)}}} \right).$$

\section{Fractal groups}
In this section, we will consider some details of the algebra $\left[\!\left[ {{P}_{\left[ 0,q \right]}} \right]\!\right]$,
$${{P}_{\left[ 0,q \right]}}=\prod\limits_{k=1}^{\infty }{\times {{P}_{0,{{q}^{k}}}}}, \qquad{{\left( P{}_{\left[ 0,q \right]} \right)}_{n,m}}=\left\{ \begin{matrix}
   1,n\left( \bmod {{q}^{k}} \right)\ge m\left( \bmod {{q}^{k}} \right),  \\
   0,n\left( \bmod {{q}^{k}} \right)<m\left( \bmod {{q}^{k}} \right).  \\
\end{matrix} \right.$$
An example of the matrix  ${{P}_{\left[ 0,q \right]}}$ is Pascal's triangle modulo 2:
$${{P}_{\left[ 0,2 \right]}}=\left(\setcounter{MaxMatrixCols}{20} \begin{matrix}
   1 & 0 & 0 & 0 & 0 & 0 & 0 & 0 & 0 & 0 & 0 & 0 & 0 & 0 & 0 & 0 & \ldots   \\
   1 & 1 & 0 & 0 & 0 & 0 & 0 & 0 & 0 & 0 & 0 & 0 & 0 & 0 & 0 & 0 & \ldots   \\
   1 & 0 & 1 & 0 & 0 & 0 & 0 & 0 & 0 & 0 & 0 & 0 & 0 & 0 & 0 & 0 & \ldots   \\
   1 & 1 & 1 & 1 & 0 & 0 & 0 & 0 & 0 & 0 & 0 & 0 & 0 & 0 & 0 & 0 & \ldots   \\
   1 & 0 & 0 & 0 & 1 & 0 & 0 & 0 & 0 & 0 & 0 & 0 & 0 & 0 & 0 & 0 & \ldots   \\
   1 & 1 & 0 & 0 & 1 & 1 & 0 & 0 & 0 & 0 & 0 & 0 & 0 & 0 & 0 & 0 & \ldots   \\
   1 & 0 & 1 & 0 & 1 & 0 & 1 & 0 & 0 & 0 & 0 & 0 & 0 & 0 & 0 & 0 & \ldots   \\
   1 & 1 & 1 & 1 & 1 & 1 & 1 & 1 & 0 & 0 & 0 & 0 & 0 & 0 & 0 & 0 & \ldots   \\
   1 & 0 & 0 & 0 & 0 & 0 & 0 & 0 & 1 & 0 & 0 & 0 & 0 & 0 & 0 & 0 & \ldots   \\
   1 & 1 & 0 & 0 & 0 & 0 & 0 & 0 & 1 & 1 & 0 & 0 & 0 & 0 & 0 & 0 & \ldots   \\
   1 & 0 & 1 & 0 & 0 & 0 & 0 & 0 & 1 & 0 & 1 & 0 & 0 & 0 & 0 & 0 & \ldots   \\
   1 & 1 & 1 & 1 & 0 & 0 & 0 & 0 & 1 & 1 & 1 & 1 & 0 & 0 & 0 & 0 & \ldots   \\
   1 & 0 & 0 & 0 & 1 & 0 & 0 & 0 & 1 & 0 & 0 & 0 & 1 & 0 & 0 & 0 & \ldots   \\
   1 & 1 & 0 & 0 & 1 & 1 & 0 & 0 & 1 & 1 & 0 & 0 & 1 & 1 & 0 & 0 & \ldots   \\
   1 & 0 & 1 & 0 & 1 & 0 & 1 & 0 & 1 & 0 & 1 & 0 & 1 & 0 & 1 & 0 & \ldots   \\
   1 & 1 & 1 & 1 & 1 & 1 & 1 & 1 & 1 & 1 & 1 & 1 & 1 & 1 & 1 & 1 & \ldots   \\
   \vdots  & \vdots  & \vdots  & \vdots  & \vdots  & \vdots  & \vdots  & \vdots  & \vdots  & \vdots  & \vdots  & \vdots  & \vdots  & \vdots  & \vdots  & \vdots  & \ddots   \\
\end{matrix} \right).$$

Denote
${{\left( {{P}_{\left[ 0,q \right]}} \right)}_{n,m}}={n\choose m}_{\left[ 0,q \right]}$. Then
$${{\left( \begin{matrix}
   {{q}^{k}}n+i  \\
   {{q}^{k}}m+j  \\
\end{matrix} \right)}_{\left[ 0,q \right]}}={{\left( \begin{matrix}
   n  \\
   m  \\
\end{matrix} \right)}_{\left[ 0,q \right]}}{{\left( \begin{matrix}
   i  \\
   j  \\
\end{matrix} \right)}_{\left[ 0,q \right]}},  \qquad 0\le i,j<{{q}^{k}},\eqno(4)$$
$${{\left( \begin{matrix}
   n  \\
   m  \\
\end{matrix} \right)}_{\left[ 0,q \right]}}=\prod\limits_{i=0}^{\infty }{{{\left( \begin{matrix}
   {{n}_{i}}  \\
   {{m}_{i}}  \\
\end{matrix} \right)}_{\left[ 0,q \right]}}},  \quad n=\sum\limits_{i=0}^{\infty }{{{n}_{i}}}{{q}^{ik}}, \quad m=\sum\limits_{i=0}^{\infty }{{{m}_{i}}{{q}^{ik}}}, \quad 0\le {{n}_{i}},{{m}_{i}}<{{q}^{k}}.$$
{\bfseries Remark 4.1.} The matrices ${{P}_{\left[ 0,q \right]}}$ were introduced in [16], where they are denoted by  ${{S}_{q}}$ and called generalized Sierpinski matrices. In [16], they are constructed as follows:
$${{S}_{q}}={{S}_{q,1}}\otimes {{S}_{q,1}}\otimes {{S}_{q,1}}\otimes ...,$$
where the sign $\otimes $ means the Kronecker product, ${{S}_{q,1}}$ is the matrix consisting of the first $q$ rows of the matrix ${{S}_{q}}$, i.e., ${{\left( {{\left( 1-x \right)}^{-1}},x \right)}_{q}}$. The property defined by identity (4) can be represented as
$${{S}_{q}}={{S}_{q,k}}\otimes {{S}_{q,k}}\otimes {{S}_{q,k}}\otimes ...,$$
where ${{S}_{q,k}}$ is the matrix consisting of the first ${{q}^{k}}$ rows of the matrix ${{S}_{q}}$. Properties of generalized Sierpinski matrices and related matrices are studied in the papers [12,14,15,16,22]. 

The algebra  $\left[\!\left[ {{P}_{\left[ 0,q \right]}} \right]\!\right]$  contains an infinite set of block groups. Since
$$\left[ {{x}^{{{q}^{k}}n+i}} \right]{{\left[ b\left( x \right) \right]}_{{{q}^{k}}}}a\left( {{x}^{{{q}^{k}}}} \right)={{a}_{n}}{{b}_{i}},  \qquad 0\le i<{{q}^{k}},$$
then 
$${{\left( {{\left[ b\left( x \right) \right]}_{{{q}^{k}}}}a\left( {{x}^{{{q}^{k}}}} \right)|{{P}_{\left[ 0,q \right]}} \right)}_{{{q}^{k}}n+i,{{q}^{k}}m+j}}=$$ 
$$={{a}_{n-m}}{{b}_{i-j}}{{\left( \begin{matrix}
   {{q}^{k}}n+i  \\
   {{q}^{k}}m+j  \\
\end{matrix} \right)}_{\left[ 0,q \right]}}={{a}_{n-m}}{{\left( \begin{matrix}
   n  \\
   m  \\
\end{matrix} \right)}_{\left[ 0,q \right]}}{{b}_{i-j}}{{\left( \begin{matrix}
   i  \\
   j  \\
\end{matrix} \right)}_{\left[ 0,q \right]}},\qquad 0\le i,j<{{q}^{k}}.$$
Thus, the matrix $\left( {{\left[ b\left( x \right) \right]}_{{{q}^{k}}}}a\left( {{x}^{{{q}^{k}}}} \right)|{{P}_{\left[ 0,q \right]}} \right)$ is a block matrix, the $\left( n,m \right)$th block of which is the matrix 
$${{a}_{n-m}}{{\left( \begin{matrix}
   n  \\
   m  \\
\end{matrix} \right)}_{\left[ 0,q \right]}}{{\left( {{\left[ b\left( x \right) \right]}_{{{q}^{k}}}}|{{P}_{\left[ 0,q \right]}} \right)}_{{{q}^{k}}}},$$
where ${{\left( {{\left[ b\left( x \right) \right]}_{{{q}^{k}}}}|{{P}_{\left[ 0,q \right]}} \right)}_{{{q}^{k}}}}$ is the matrix consisting of the first ${{q}^{k}}$ rows of the matrix \linebreak$\left( {{\left[ b\left( x \right) \right]}_{{{q}^{k}}}}|{{P}_{\left[ 0,q \right]}} \right)$. For example,
$$\left( {{\left[ b\left( x \right) \right]}_{2}}a\left( {{x}^{2}} \right)|{{P}_{\left[ 0,2 \right]}} \right)=\left( \begin{matrix}
   {{a}_{0}}{{b}_{0}} & 0 & 0 & 0 & 0 & 0 & 0 & 0 & \ldots   \\
   {{a}_{0}}{{b}_{1}} & {{a}_{0}}{{b}_{0}} & 0 & 0 & 0 & 0 & 0 & 0 & \ldots   \\
   {{a}_{1}}{{b}_{0}} & 0 & {{a}_{0}}{{b}_{0}} & 0 & 0 & 0 & 0 & 0 & \ldots   \\
   {{a}_{1}}{{b}_{1}} & {{a}_{1}}{{b}_{0}} & {{a}_{0}}{{b}_{1}} & {{a}_{0}}{{b}_{0}} & 0 & 0 & 0 & 0 & \ldots   \\
   {{a}_{2}}{{b}_{0}} & 0 & 0 & 0 & {{a}_{0}}{{b}_{0}} & 0 & 0 & 0 & \ldots   \\
   {{a}_{2}}{{b}_{1}} & {{a}_{2}}{{b}_{0}} & 0 & 0 & {{a}_{0}}{{b}_{1}} & {{a}_{0}}{{b}_{0}} & 0 & 0 & \ldots   \\
   {{a}_{3}}{{b}_{0}} & 0 & {{a}_{2}}{{b}_{0}} & 0 & {{a}_{1}}{{b}_{0}} & 0 & {{a}_{0}}{{b}_{0}} & 0 & \ldots   \\
   {{a}_{3}}{{b}_{1}} & {{a}_{3}}{{b}_{0}} & {{a}_{2}}{{b}_{1}} & {{a}_{2}}{{b}_{0}} & {{a}_{1}}{{b}_{1}} & {{a}_{1}}{{b}_{0}} & {{a}_{0}}{{b}_{1}} & {{a}_{0}}{{b}_{0}} & \ldots   \\
   \vdots  & \vdots  & \vdots  & \vdots  & \vdots  & \vdots  & \vdots  & \vdots  & \ddots   \\
\end{matrix} \right),$$
$$\left( {{\left[ b\left( x \right) \right]}_{4}}a\left( {{x}^{4}} \right)|{{P}_{\left[ 0,2 \right]}} \right)=\left( \begin{matrix}
   {{a}_{0}}{{b}_{0}} & 0 & 0 & 0 & 0 & 0 & 0 & 0 & \ldots   \\
   {{a}_{0}}{{b}_{1}} & {{a}_{0}}{{b}_{0}} & 0 & 0 & 0 & 0 & 0 & 0 & \ldots   \\
   {{a}_{0}}{{b}_{2}} & 0 & {{a}_{0}}{{b}_{0}} & 0 & 0 & 0 & 0 & 0 & \ldots   \\
   {{a}_{0}}{{b}_{3}} & {{a}_{0}}{{b}_{2}} & {{a}_{0}}{{b}_{1}} & {{a}_{0}}{{b}_{0}} & 0 & 0 & 0 & 0 & \ldots   \\
   {{a}_{1}}{{b}_{0}} & 0 & 0 & 0 & {{a}_{0}}{{b}_{0}} & 0 & 0 & 0 & \ldots   \\
   {{a}_{1}}{{b}_{1}} & {{a}_{1}}{{b}_{0}} & 0 & 0 & {{a}_{0}}{{b}_{1}} & {{a}_{0}}{{b}_{0}} & 0 & 0 & \ldots   \\
   {{a}_{1}}{{b}_{2}} & 0 & {{a}_{1}}{{b}_{0}} & 0 & {{a}_{0}}{{b}_{2}} & 0 & {{a}_{0}}{{b}_{0}} & 0 & \ldots   \\
   {{a}_{1}}{{b}_{3}} & {{a}_{1}}{{b}_{2}} & {{a}_{1}}{{b}_{1}} & {{a}_{1}}{{b}_{0}} & {{a}_{0}}{{b}_{3}} & {{a}_{0}}{{b}_{2}} & {{a}_{0}}{{b}_{1}} & {{a}_{0}}{{b}_{0}} & \ldots   \\
   \vdots  & \vdots  & \vdots  & \vdots  & \vdots  & \vdots  & \vdots  & \vdots  & \ddots   \\
\end{matrix} \right).$$
Hence,  if ${{a}_{0}}\ne 0$, ${{b}_{0}}\ne 0$, these matrices form a group whose elements are multiplied by the rule
$$\left( {{\left[ {{b}_{1}}\left( x \right) \right]}_{{{q}^{k}}}}{{a}_{1}}\left( {{x}^{{{q}^{k}}}} \right)|{{P}_{0,q}} \right)\left( {{\left[ {{b}_{2}}\left( x \right) \right]}_{{{q}^{k}}}}{{a}_{2}}\left( {{x}^{{{q}^{k}}}} \right)|{{P}_{0,q}} \right)=$$
$$=\left( {{\left[ {{b}_{1}}\left( x \right)\circ {{b}_{2}}\left( x \right) \right]}_{{{q}^{k}}}}\left( {{a}_{1}}\left( {{x}^{{{q}^{k}}}} \right)\circ {{a}_{2}}\left( {{x}^{{{q}^{k}}}} \right) \right)|{{P}_{0,q}} \right),$$
${{\left[ {{b}_{1}}\left( x \right)\circ {{b}_{2}}\left( x \right) \right]}_{q}}={{\left[ {{b}_{1}}\left( x \right){{b}_{2}}\left( x \right) \right]}_{q}}$. We will denote this group by ${{B}_{q,k}}$. Let's turn to a family of the "fractal" series ${}^{q}a\left( x \right)$, such that
$${}^{q}a\left( x \right)={{\left[ {}^{q}a\left( x \right) \right]}_{q}}{}^{q}a\left( {{x}^{q}} \right)={{\left[ {}^{q}a\left( x \right) \right]}_{{{q}^{k}}}}{}^{q}a\left( {{x}^{{{q}^{k}}}} \right)=\prod\limits_{m=0}^{\infty }{\left( \sum\limits_{n=0}^{{{q}^{k}}-1}{{}^{q}{{a}_{n}}{{x}^{n{{q}^{mk}}}}} \right)},$$
$${}^{q}{{a}_{0}}=1,  \qquad{}^{q}{{a}_{{{q}^{k}}n+i}}={}^{q}{{a}_{n}}{}^{q}{{a}_{i}},  \qquad 0\le i<{{q}^{k}},$$ 
$${}^{q}{{a}_{n}}=\prod\limits_{i=0}^{\infty }{{}^{q}{{a}_{{{n}_{i}}}}}, \qquad n=\sum\limits_{i=0}^{\infty }{{{n}_{i}}}{{q}^{ik}}={{n}_{0}}+{{q}^{k}}\left( {{n}_{1}}+{{q}^{k}}\left( {{n}_{2}}+... \right) \right),  \qquad 0\le {{n}_{i}}<{{q}^{k}}.$$
Considering the properties of the coefficients  ${n\choose m}_{\left[ 0,q \right]}$, we have:
$${{\left( {}^{q}a\left( x \right)|P{}_{\left[ 0,q \right]} \right)}_{{{q}^{k}}n+i,{{q}^{k}}m+j}}={{\left( {}^{q}a\left( x \right)|P{}_{\left[ 0,q \right]} \right)}_{n,m}}{{\left( {}^{q}a\left( x \right)|P{}_{\left[ 0,q \right]} \right)}_{i,j}},\qquad 0\le i,j<{{q}^{k}},$$
$${{\left( {}^{q}a\left( x \right)|P{}_{\left[ 0,q \right]} \right)}_{n,m}}=\prod\limits_{i=0}^{\infty }{{{\left( {}^{q}a\left( x \right)|P{}_{\left[ 0,q \right]} \right)}_{{{n}_{i}},{{m}_{i}}}}},\, 
 n=\sum\limits_{i=0}^{\infty }{{{n}_{i}}}{{q}^{ik}},\, m=\sum\limits_{i=0}^{\infty }{{{m}_{i}}}{{q}^{ik}},\, 0\le {{n}_{i}},{{m}_{i}}<{{q}^{k}}.$$
Denote $\left[ n,\to  \right]\left( {}^{q}a\left( x \right)|P{}_{\left[ 0,q \right]} \right)={{u}_{n}}\left( x \right)$. Then
$${{u}_{n}}\left( x \right)=\sum\limits_{m=0}^{n}{{}^{q}{{a}_{n-m}}{{x}^{m}}}, \quad 0\le n<q;  \qquad{{u}_{{{q}^{k}}n+i}}\left( x \right)={{u}_{n}}\left( {{x}^{{{q}^{k}}}} \right){{u}_{i}}\left( x \right), \quad 0\le i<{{q}^{k}};$$
$${{u}_{n}}\left( x \right)=\prod\limits_{i=0}^{\infty }{{{u}_{{{n}_{i}}}}}\left( {{x}^{{{q}^{ik}}}} \right),  \qquad n=\sum\limits_{i=0}^{\infty }{{{n}_{i}}}{{q}^{ik}},  \qquad 0\le {{n}_{i}}<{{q}^{k}}.$$
{\bfseries Theorem 4.1.} \emph{ The subalgebra of the matrices $\left( a\left( {{x}^{{{q}^{k}}}} \right)|{{P}_{\left[ 0,q \right]}} \right)$ is isomorphic to the algebra $\left[\!\left[ {{P}_{\left[ 0,q \right]}} \right]\!\right]$. }\\
{\bfseries Proof.} Rows of the matrix ${{P}_{\left[ 0,q \right]}}$  are related by the relation $\left[ n,\to  \right]{{P}_{\left[ 0,q \right]}}={{u}_{n}}\left( x \right)$, $\left[ {{q}^{k}}n,\to  \right]{{P}_{\left[ 0,q \right]}}={{u}_{n}}\left( {{x}^{{{q}^{k}}}} \right)$. Then
$$\left[ n,\to  \right]\left( a\left( x \right)|{{P}_{\left[ 0,q \right]}} \right)={{a}_{n}}\left( x \right),   \qquad\left[ {{q}^{k}}n,\to  \right]\left( a\left( {{x}^{{{q}^{k}}}} \right)|{{P}_{\left[ 0,q \right]}} \right)={{a}_{n}}\left( {{x}^{{{q}^{k}}}} \right),$$
$$\left[ {{x}^{n}} \right]a\left( x \right)\circ b\left( x \right)=\left[ {{x}^{{{q}^{k}}n}} \right]a\left( {{x}^{{{q}^{k}}}} \right)\circ b\left( {{x}^{{{q}^{k}}}} \right). \qquad\square $$
{\bfseries Theorem 4.2.} \emph{The matrices $\left( {}^{q}a\left( x \right)|{{P}_{\left[ 0,q \right]}} \right)$ form a group (denoted by ${{F}_{q}}$) isomorphic to the group of the matrices ${{\left( {}^{q}a\left( x \right),x \right)}_{q}}$.}\\
{\bfseries Proof. } The matrix $\left( {}^{q}a\left( x \right)|{{P}_{\left[ 0,q \right]}} \right)$ is an element of each group ${{B}_{q,k}}$. Hence, taking into account Theorem 4.1., 
$${}^{q}a\left( x \right)\circ {}^{q}b\left( x \right)={{\left[ {}^{q}a\left( x \right)\circ {}^{q}b\left( x \right) \right]}_{{{q}^{k}}}}\left( {}^{q}a\left( {{x}^{{{q}^{k}}}} \right)\circ {}^{q}b\left( {{x}^{{{q}^{k}}}} \right) \right)={}^{q}c\left( x \right).$$
It remains to add that ${{\left( {}^{q}a\left( x \right)|{{P}_{\left[ 0,q \right]}} \right)}_{q}}={{\left( {}^{q}a\left( x \right),x \right)}_{q}}$.\qquad  $\square $\\ 
{\bfseries Theorem 4.3.} \emph{The matrices $\left( {}^{q}a\left( {{x}^{{{q}^{k}}}} \right)|{{P}_{\left[ 0,{{q}^{k+1}} \right]}} \right)$ form a subgroup in ${{F}_{{{q}^{k+1}}}}$ isomorphic to the group ${{F}_{q}}$.}\\
{\bfseries  Proof.} Since 
$${}^{q}a\left( {{x}^{{{q}^{k}}}} \right)=\left( \sum\limits_{n=0}^{q-1}{{}^{q}{{a}_{n}}{{x}^{n{{q}^{k}}}}} \right){}^{q}a\left( {{x}^{{{q}^{k+1}}}} \right)={{\left[ {}^{q}a\left( {{x}^{{{q}^{k}}}} \right) \right]}_{{{q}^{k+1}}}}{}^{q}a\left( {{x}^{{{q}^{k+1}}}} \right),$$
then $\left( {}^{q}a\left( {{x}^{{{q}^{k}}}} \right)|{{P}_{\left[ 0,{{q}^{k+1}} \right]}} \right)\in {{F}_{{{q}^{k+1}}}}$. It remains to note that the group of the matrices ${{\left( {}^{q}a\left( {{x}^{{{q}^{k}}}} \right),x \right)}_{{{q}^{k+1}}}}$ is isomorphic to the group of the matrices ${{\left( {}^{q}a\left( x \right),x \right)}_{q}}$.\qquad $\square $\\ 
{\bfseries Theorem 4.4. } \emph{$$\log \circ {}^{q}a\left( x \right)=\sum\limits_{k=0}^{\infty }{{{\left[ \log {}^{q}a\left( {{x}^{{{q}^{k}}}} \right) \right]}_{{{q}^{k+1}}}}=}\sum\limits_{k=0}^{\infty }{\sum\limits_{n=1}^{q-1}{{{l}_{n}}{{x}^{n{{q}^{k}}}}}}=\sum\limits_{n=1}^{q-1}{{{l}_{n}}\sum\limits_{k=0}^{\infty }{{{x}^{n{{q}^{k}}}}}},$$
where ${{l}_{n}}=\left[ {{x}^{n}} \right]\log {}^{q}a\left( x \right)$.}\\
{\bfseries Proof.} Denote $\log \circ {}^{q}a\left( x \right)={{l}_{q}}\left( x \right)$. Since $\log \circ {{\left[ {}^{q}a\left( x \right) \right]}_{q}}={{\left[ \log {}^{q}a\left( x \right) \right]}_{q}}$ and by \linebreak Theorem 4.1. $\log \circ {}^{q}a\left( {{x}^{{{q}^{k}}}} \right)={{l}_{q}}\left( {{x}^{{{q}^{k}}}} \right)$, then
$${{l}_{q}}\left( x \right)={{\left[ \log {}^{q}a\left( x \right) \right]}_{q}}+{{l}_{q}}\left( {{x}^{q}} \right)=\sum\limits_{k=0}^{\infty }{{{\left[ \log {}^{q}a\left( {{x}^{{{q}^{k}}}} \right) \right]}_{{{q}^{k+1}}}}}.  \qquad\square $$

Thus,
$$\log P{}_{\left[ 0,q \right]}=\left( \sum\limits_{k=0}^{\infty }{\sum\limits_{n=1}^{q-1}{\frac{{{x}^{n{{q}^{k}}}}}{n}}}|P{}_{\left[ 0,q \right]} \right).$$
\section{The group $R\left( {{P}_{0}} \right)$}
Denote
$\left( a\left( x \right)|P{}_{0} \right)={{\left( a\left( x \right),1 \right)}_{0}}$,
where the form of the matrix  ${{P}_{0}}$ is indicated separately. Let's construct the matrix ${{\left( 1,a\left( x \right) \right)}_{0}}$ by the rule
$${{\left( 1,a\left( x \right) \right)}_{0}}{{x}^{n}}={{\left( {{a}^{\left( n \right)}}\left( x \right),1 \right)}_{0}}{{x}^{n}}={{x}^{n}}\circ {{a}^{\left( n \right)}}\left( x \right).$$
Denote
$${{\left( 1,a\left( x \right) \right)}_{0}}b\left( x \right)=b{{\left( a\left( x \right) \right)}_{0}},\qquad
{{\left( b\left( x \right),1 \right)}_{0}}{{\left( 1,a\left( x \right) \right)}_{0}}={{\left( b\left( x \right),a\left( x \right) \right)}_{0}}.$$
{\bfseries Theorem  5.1.} \emph{ The matrices ${{\left( b\left( x \right),a\left( x \right) \right)}_{0}}$, ${{b}_{0}}\ne 0$, ${{a}_{0}}\ne 0$ form a group whose elements are multiplied by the rule}
$${{\left( b\left( x \right),a\left( x \right) \right)}_{0}}{{\left( f\left( x \right),g\left( x \right) \right)}_{0}}={{\left( b\left( x \right)\circ f{{\left( a\left( x \right) \right)}_{0}},a\left( x \right)\circ g{{\left( a\left( x \right) \right)}_{0}} \right)}_{0}}.$$
{\bfseries Proof.} Since
$${{x}^{m}}\circ {{x}^{n}}={{\left( \begin{matrix}
   m+n  \\
   n  \\
\end{matrix} \right)}_{0}}{{x}^{m+n}} ,  \qquad{{x}^{m}}\circ b\left( x \right)=\sum\limits_{n=0}^{\infty }{{{b}_{n}}}{{\left( \begin{matrix}
   m+n  \\
   n  \\
\end{matrix} \right)}_{0}}{{x}^{m+n}},$$
$${{\left( 1,a\left( x \right) \right)}_{0}}\left( {{x}^{m}}\circ b\left( x \right) \right)=\sum\limits_{n=0}^{\infty }{{{b}_{n}}}{{\left( \begin{matrix}
   m+n  \\
   n  \\
\end{matrix} \right)}_{0}}{{x}^{m+n}}\circ {{a}^{\left( m+n \right)}}\left( x \right)=$$
$$={{x}^{m}}\circ {{a}^{\left( m \right)}}\left( x \right)\circ \sum\limits_{n=0}^{\infty }{{{b}_{n}}}{{x}^{n}}\circ {{a}^{\left( n \right)}}\left( x \right)={{x}^{m}}\circ {{a}^{\left( m \right)}}\left( x \right)\circ b{{\left( a\left( x \right) \right)}_{0}},$$
then
$${{\left( 1,a\left( x \right) \right)}_{0}}{{\left( b\left( x \right),1 \right)}_{0}}={{\left( b{{\left( a\left( x \right) \right)}_{0}},a\left( x \right) \right)}_{0}}.$$
Then
$${{\left( 1,a\left( x \right) \right)}_{0}}\left( b\left( x \right)\circ c\left( x \right) \right)=b{{\left( a\left( x \right) \right)}_{0}}\circ c{{\left( a\left( x \right) \right)}_{0}},$$
$${{\left( 1,a\left( x \right) \right)}_{0}}\left( {{x}^{m}}\circ {{b}^{\left( m \right)}}\left( x \right) \right)={{x}^{m}}\circ {{a}^{\left( m \right)}}\left( x \right)\circ {{\left( b{{\left( a\left( x \right) \right)}_{0}} \right)}^{\left( m \right)}},$$
or
$${{\left( 1,a\left( x \right) \right)}_{0}}{{\left( 1,b\left( x \right) \right)}_{0}}={{\left( 1,a\left( x \right)\circ b{{\left( a\left( x \right) \right)}_{0}} \right)}_{0}}.  \qquad\square $$ 

 The generalized Riordan group with the central element ${{P}_{c\left( x \right)}}$  will be denoted by $R\left( {{P}_{c\left( x \right)}} \right)$; the   group of the matrices ${{\left( b\left( x \right),a\left( x \right) \right)}_{0}}$ will be denoted by $R\left( {{P}_{0}} \right)$. The terminology used for the Riordan group will be used for the group $R\left( {{P}_{0}} \right)$. The subgroups of the matrices ${{\left( b\left( x \right),1 \right)}_{0}}$, ${{\left( 1,a\left( x \right) \right)}_{0}}$ will be called the Appell subgroup and the Lagrange subgroup, respectively. The subgroup of the matrices ${{\left( a\left( x \right),a\left( x \right) \right)}_{0}}$ isomorphic to the Lagrange subgroup will be called the Bell subgroup. Since $x\circ a\left( x \right)=x\circ b\left( x \right)$ for an infinite set of series $b\left( x \right)$, then the matrices ${{\left( a\left( x \right),a\left( x \right) \right)}_{0}}$ have an advantage over the matrices ${{\left( 1,a\left( x \right) \right)}_{0}}$ in practical terms. The following example illustrates this.\\
 {\bfseries Example 5.1.} Here ${{P}_{0}}={{P}_{0,2}}$;
$${{\left( 1,\frac{1}{1-{{x}^{2}}} \right)}_{0}}{{\left( 1,\frac{1}{1-x} \right)}_{0}}={{\left( 1,\frac{1}{1-x-{{x}^{2}}} \right)}_{0}},$$
$$\left( \begin{matrix}
   1 & 0 & 0 & 0 & 0 & 0 & 0 & \cdots   \\
   0 & 1 & 0 & 0 & 0 & 0 & 0 & \cdots   \\
   0 & 0 & 1 & 0 & 0 & 0 & 0 & \cdots   \\
   0 & 1 & 0 & 1 & 0 & 0 & 0 & \cdots   \\
   0 & 0 & 2 & 0 & 1 & 0 & 0 & \cdots   \\
   0 & 1 & 0 & 3 & 0 & 1 & 0 & \cdots   \\
   0 & 0 & 3 & 0 & 4 & 0 & 1 & \cdots   \\
   \vdots  & \vdots  & \vdots  & \vdots  & \vdots  & \vdots  & \vdots  & \ddots   \\
\end{matrix} \right)\!\left( \begin{matrix}
   1 & 0 & 0 & 0 & 0 & 0 & 0 & \cdots   \\
   0 & 1 & 0 & 0 & 0 & 0 & 0 & \cdots   \\
   0 & 0 & 1 & 0 & 0 & 0 & 0 & \cdots   \\
   0 & 1 & 2 & 1 & 0 & 0 & 0 & \cdots   \\
   0 & 0 & 2 & 0 & 1 & 0 & 0 & \cdots   \\
   0 & 1 & 4 & 3 & 4 & 1 & 0 & \cdots   \\
   0 & 0 & 3 & 0 & 4 & 0 & 1 & \cdots   \\
   \vdots  & \vdots  & \vdots  & \vdots  & \vdots  & \vdots  & \vdots  & \ddots   \\
\end{matrix} \right)\!=\!\left( \begin{matrix}
   1 & 0 & 0 & 0 & 0 & 0 & 0 & \cdots   \\
   0 & 1 & 0 & 0 & 0 & 0 & 0 & \cdots   \\
   0 & 0 & 1 & 0 & 0 & 0 & 0 & \cdots   \\
   0 & 2 & 2 & 1 & 0 & 0 & 0 & \cdots   \\
   0 & 0 & 4 & 0 & 1 & 0 & 0 & \cdots   \\
   0 & 5 & 10 & 6 & 4 & 1 & 0 & \cdots   \\
   0 & 0 & 14 & 0 & 8 & 0 & 1 & \cdots   \\
   \vdots  & \vdots  & \vdots  & \vdots  & \vdots  & \vdots  & \vdots  & \ddots   \\
\end{matrix} \right);$$
$${{\left( \frac{1}{1-{{x}^{2}}},\frac{1}{1-{{x}^{2}}} \right)}_{0}}{{\left( \frac{1}{1-x},\frac{1}{1-x} \right)}_{0}}={{\left( \frac{1}{1-x-{{x}^{2}}},\frac{1}{1-x-{{x}^{2}}} \right)}_{0}},$$
$$\left( \begin{matrix}
   1 & 0 & 0 & 0 & 0 & 0 & \cdots   \\
   0 & 1 & 0 & 0 & 0 & 0 & \cdots   \\
   1 & 0 & 1 & 0 & 0 & 0 & \cdots   \\
   0 & 2 & 0 & 1 & 0 & 0 & \cdots   \\
   1 & 0 & 3 & 0 & 1 & 0 & \cdots   \\
   0 & 3 & 0 & 4 & 0 & 1 & \cdots   \\
   \vdots  & \vdots  & \vdots  & \vdots  & \vdots  & \vdots  & \ddots   \\
\end{matrix} \right)\left( \begin{matrix}
   1 & 0 & 0 & 0 & 0 & 0 & \cdots   \\
   1 & 1 & 0 & 0 & 0 & 0 & \cdots   \\
   1 & 0 & 1 & 0 & 0 & 0 & \cdots   \\
   1 & 2 & 3 & 1 & 0 & 0 & \cdots   \\
   1 & 0 & 3 & 0 & 1 & 0 & \cdots   \\
   1 & 3 & 9 & 4 & 5 & 1 & \cdots   \\
     \ & \vdots  & \vdots  & \vdots  & \vdots  & \vdots  & \vdots  & \ddots   \\
\end{matrix} \right)=\left( \begin{matrix}
   1 & 0 & 0 & 0 & 0 & 0 & \cdots   \\
   1 & 1 & 0 & 0 & 0 & 0 & \cdots   \\
   2 & 0 & 1 & 0 & 0 & 0 & \cdots   \\
   3 & 4 & 3 & 1 & 0 & 0 & \cdots   \\
   5 & 0 & 6 & 0 & 1 & 0 & \cdots   \\
   8 & 14 & 21 & 8 & 5 & 1 & \cdots   \\
   \vdots  & \vdots  & \vdots  & \vdots  & \vdots  & \vdots  & \ddots   \\
\end{matrix} \right).$$

This example is also an illustration of a intersection of the group $R\left( {{P}_{0,q}} \right)$ with the ordinary Riordan group:  ${{\left( b\left( {{x}^{q}} \right),a\left( {{x}^{q}} \right) \right)}_{0}}=\left( b\left( {{x}^{q}} \right),xa\left( {{x}^{q}} \right) \right)$.

 Note thet $R\left( {{P}_{c\left( x \right)}} \right)\bigcap R\left( {{P}_{0}} \right)={{\left( 1,\varphi  \right)}_{0}}=\left( 1,\varphi x \right)$ for any ${{P}_{c\left( x \right)}}$ and ${{P}_{0}}$. Thus, $a{{\left( \varphi  \right)}_{0}}=a\left( \varphi x \right)$,
$${{\left( 1,\varphi  \right)}_{0}}{{\left( b\left( x \right),a\left( x \right) \right)}_{0}}={{\left( b\left( \varphi x \right),\varphi a\left( \varphi x \right) \right)}_{0}},
\qquad{{\left( b\left( x \right),a\left( x \right) \right)}_{0}}{{\left( 1,\varphi  \right)}_{0}}={{\left( b\left( x \right),\varphi a\left( x \right) \right)}_{0}}.$$

The following theorem is an analog of the Lagrange inversion theorem in its generalized form.\\
{\bfseries Theorem 5.2.} \emph{Each formal power series $a\left( x \right)\in \left[\!\left[ {{P}_{0}},a\left( x \right) \right]\!\right]$, ${{a}_{0}}=1$,  is associated with a family of series $_{\left( \beta  \right)}a\left( x \right)\in \left[\!\left[ {{P}_{0}},a\left( x \right) \right]\!\right]$, ${}_{\left( 0 \right)}a\left( x \right)=a\left( x \right)$, such that}
$${}_{\left( \beta  \right)}a{{\left( {{a}^{\left( -\beta  \right)}}\left( x \right) \right)}_{0}}=a\left( x \right),   \qquad a{{\left( {}_{\left( \beta  \right)}{{a}^{\left( \beta  \right)}}\left( x \right) \right)}_{0}}={}_{\left( \beta  \right)}a\left( x \right),$$
$$\left[ {{x}^{n}} \right]{}_{\left( \beta  \right)}{{a}^{\left( \varphi  \right)}}\left( x \right)=\left[ {{x}^{n}} \right]\left( 1-x\beta {{\left( \log \circ a\left( x \right) \right)}^{\prime }} \right)\circ {{a}^{\left( \varphi +\beta n \right)}}\left( x \right)
=\frac{\varphi }{\varphi +\beta n}\left[ {{x}^{n}} \right]{{a}^{\left( \varphi +\beta n \right)}}\left( x \right),$$
$$\left[ {{x}^{n}} \right]\left( 1+x\beta {{\left( \log \circ {}_{\left( \beta  \right)}a\left( x \right) \right)}^{\prime }} \right)\circ {}_{\left( \beta  \right)}{{a}^{\left( \varphi  \right)}}\left( x \right)
=\frac{\varphi +\beta n}{\varphi }\left[ {{x}^{n}} \right]{}_{\left( \beta  \right)}{{a}^{\left( \varphi  \right)}}\left( x \right)=\left[ {{x}^{n}} \right]{{a}^{\left( \varphi +\beta n \right)}}\left( x \right).$$
{\bfseries Proof.} If the matrices ${{\left( 1,{{a}^{\left( -1 \right)}}\left( x \right) \right)}_{0}}$, ${{a}_{0}}=1$, ${{\left( 1,b\left( x \right) \right)}_{0}}$, ${{b}_{0}}=1$, are mutually inverse, then 
$${{\left( 1,{{a}^{\left( -1 \right)}}\left( x \right) \right)}_{0}}b\left( x \right)=a\left( x \right),   \qquad{{\left( 1,b\left( x \right) \right)}_{0}}a\left( x \right)=b\left( x \right).$$
Since
$$x{{\left( {{x}^{n}}\circ {{b}^{\left( n \right)}}\left( x \right) \right)}^{\prime }}={{x}^{n}}\circ x{{\left( {{b}^{\left( n \right)}}\left( x \right) \right)}^{\prime }}+x{{\left( {{x}^{n}} \right)}^{\prime }}\circ {{b}^{\left( n \right)}}\left( x \right)=$$
$$={{x}^{n}}\circ n{{b}^{\left( n-1 \right)}}\left( x \right)\circ x{{b}^{\prime }}\left( x \right)+n{{x}^{n}}\circ {{b}^{\left( n \right)}}\left( x \right)
=n{{x}^{n}}\circ {{b}^{\left( n \right)}}\left( x \right)\circ \left( 1+x{{\left( \log \circ b\left( x \right) \right)}^{\prime }} \right),$$
then
$$\left( x,x \right)D{{\left( 1,b\left( x \right) \right)}_{0}}={{\left( \left( 1+x{{\left( \log \circ b\left( x \right) \right)}^{\prime }} \right),b\left( x \right) \right)}_{0}}\left( x,x \right)D,$$
$${{\left( 1,b\left( x \right) \right)}_{0}}x{a}'\left( x \right)={{\left( 1+x{{\left( \log \circ b\left( x \right) \right)}^{\prime }} \right)}^{\left( -1 \right)}}\circ x{b}'\left( x \right).$$
From here we find:
$$\left( 1+x{{\left( \log \circ b\left( x \right) \right)}^{\prime }},b\left( x \right) \right)_{0}^{-1}={{\left( 1-x{{\left( \log \circ a\left( x \right) \right)}^{\prime }},{{a}^{\left( -1 \right)}}\left( x \right) \right)}_{0}}.$$
Denote
$$\left[ {{x}^{n}} \right]{{a}^{\left( m \right)}}\left( x \right)=a_{n}^{\left( m \right)},  \qquad\left[ {{x}^{n}} \right]\left( 1-x{{\left( \log \circ a\left( x \right) \right)}^{\prime }} \right)\circ {{a}^{\left( m \right)}}\left( x \right)=c_{n}^{\left( m \right)},$$ 
 $${{a}_{m}}\left( x \right)=\sum\limits_{n=0}^{\infty }{a_{n}^{\left( m+n \right)}}{{x}^{n}},  \qquad{{c}_{m}}\left( x \right)=\sum\limits_{n=0}^{\infty }{c_{n}^{\left( m+n \right)}{{x}^{n}}}.$$
We construct the matrix $A$ whose $m$th column has the generating function  ${{x}^{m}}\circ {{a}_{m}}\left( x \right)$ and the matrix  $C$ whose $m$th column has the generating function ${{x}^{m}}\circ {{c}_{m}}\left( x \right)$:
$$A=\left( \begin{matrix}
   a_{0}^{\left( 0 \right)} & 0 & 0 & 0 & \cdots   \\
   a_{1}^{\left( 1 \right)} & a_{0}^{\left( 1 \right)} & 0 & 0 & \cdots   \\
   a_{2}^{\left( 2 \right)} & a_{1}^{\left( 2 \right)} & a_{0}^{\left( 2 \right)} & 0 & \cdots   \\
   a_{3}^{\left( 3 \right)} & a_{2}^{\left( 3 \right)} & a_{1}^{\left( 3 \right)} & a_{0}^{\left( 3 \right)} & \cdots   \\
   \vdots  & \vdots  & \vdots  & \vdots  & \ddots   \\
\end{matrix} \right)\times {{P}_{0}},  \qquad C=\left( \begin{matrix}
   c_{0}^{\left( 0 \right)} & 0 & 0 & 0 & \cdots   \\
   c_{1}^{\left( 1 \right)} & c_{0}^{\left( 1 \right)} & 0 & 0 & \cdots   \\
   c_{2}^{\left( 2 \right)} & c_{1}^{\left( 2 \right)} & c_{0}^{\left( 2 \right)} & 0 & \cdots   \\
   c_{3}^{\left( 3 \right)} & c_{2}^{\left( 3 \right)} & c_{1}^{\left( 3 \right)} & c_{0}^{\left( 3 \right)} & \cdots   \\
   \vdots  & \vdots  & \vdots  & \vdots  & \ddots   \\
\end{matrix} \right)\times {{P}_{0}}.$$
It's obvious that
$$\left[ n,\to  \right]A=\left[ n,\to  \right]{{\left( {{a}^{\left( n \right)}}\left( x \right),1 \right)}_{0}},\quad
\left[ n,\to  \right]C=\left[ n,\to  \right]{{\left( \left( 1-x{{\left( \log a\left( x \right) \right)}^{\prime }} \right)\circ {{a}^{\left( n \right)}}\left( x \right),1 \right)}_{0}}.$$
Since
$$\left( 1-x{a}'\left( x \right)\circ {{a}^{\left( -1 \right)}}\left( x \right) \right)\circ {{a}^{\left( m \right)}}\left( x \right)={{a}^{\left( m \right)}}\left( x \right)-\frac{x}{m}{{\left( {{a}^{\left( m \right)}}\left( x \right) \right)}^{\prime }},$$
or
$$\left[ {{x}^{n}} \right]\left( 1-x{{\left( \log \circ a\left( x \right) \right)}^{\prime }} \right)\circ {{a}^{\left( m \right)}}\left( x \right)=\frac{m-n}{m}\left[ {{x}^{n}} \right]{{a}^{\left( m \right)}}\left( x \right),$$
then
$$\left[ {{x}^{m+n}} \right]A\left( {{x}^{m}}\circ \left( 1-x{{\left( \log \circ a\left( x \right) \right)}^{\prime }} \right)\circ {{a}^{\left( -m \right)}}\left( x \right) \right)=\left[ {{x}^{m+n}} \right]C\left( {{x}^{m}}\circ {{a}^{\left( -m \right)}}\left( x \right) \right)=$$
$$=\left[ {{x}^{n}} \right]\left( 1-x{{\left( \log \circ a\left( x \right) \right)}^{\prime }} \right)\circ {{a}^{\left( n \right)}}\left( x \right)=\left\{ \begin{matrix}
   1,n=0,  \\
   0,n>0.  \\
\end{matrix} \right.$$
Thus,
$$A={{\left( 1+x{{\left( \log \circ b\left( x \right) \right)}^{\prime }},b\left( x \right) \right)}_{0}},   \qquad C={{\left( 1,b\left( x \right) \right)}_{0}},$$ 
$$\left[ {{x}^{n}} \right]\left( 1+x{{\left( \log \circ b\left( x \right) \right)}^{\prime }} \right)\circ {{b}^{\left( m \right)}}\left( x \right)=\frac{m+n}{m}\left[ {{x}^{n}} \right]{{b}^{\left( m \right)}}\left( x \right)=\left[ {{x}^{n}} \right]{{a}^{\left( m+n \right)}}\left( x \right),$$
 $$\left[ {{x}^{n}} \right]{{b}^{\left( m \right)}}\left( x \right)=\left[ {{x}^{n}} \right]\left( 1-x{{\left( \log \circ a\left( x \right) \right)}^{\prime }} \right)\circ {{a}^{\left( m+n \right)}}\left( x \right)=\frac{m}{m+n}\left[ {{x}^{n}} \right]{{a}^{\left( m+n \right)}}\left( x \right).$$
Denote
$$\left( 1,{{a}^{\left( -\beta  \right)}}\left( x \right) \right)_{0}^{-1}={{\left( 1,{}_{\left( \beta  \right)}{{a}^{\left( \beta  \right)}}\left( x \right) \right)}_{0}}.$$
Then
$$\left[ {{x}^{n}} \right]{}_{\left( \beta  \right)}{{a}^{\left( \beta m \right)}}\left( x \right)=\frac{\beta m}{\beta m+\beta n}\left[ {{x}^{n}} \right]{{a}^{\left( \beta m+\beta n \right)}}\left( x \right).$$
Let ${{c}_{n}}\left( \varphi  \right)$ be the convolution polynomials of the series $a\left( x \right)$.  Then
$$_{\left( \beta  \right)}{{a}^{\left( \varphi  \right)}}\left( x \right)=\sum\limits_{n=0}^{\infty }{\frac{\varphi }{\varphi +\beta n}}{{c}_{n}}\left( \varphi +\beta n \right){{x}^{n}}. \qquad\square $$ 
{\bfseries Example 5.2. } Consider an analog of the exponential series from the family of series $^{q}a\left( x \right)$:
$$^{q}\varepsilon \left( x \right)={{\left[ {{e}^{x}} \right]}_{q}}{}^{q}\varepsilon \left( {{x}^{q}} \right) ,  \qquad^{q}\varepsilon \left( x \right)\in \left[\!\left[ {{P}_{\left[ 0,q \right]}},a\left( x \right) \right]\!\right],  \qquad\log \circ {}^{q}\varepsilon \left( x \right)=\sum\limits_{k=0}^{\infty }{{{x}^{{{q}^{k}}}}},$$
$$\left[ {{x}^{n}} \right]{}^{q}{{\varepsilon }^{\left( \varphi  \right)}}\left( x \right)=\frac{{{\varphi }^{\left\{ n \right\}}}}{\left( n \right)!} , \quad\left\{ n \right\}=\sum\limits_{i=0}^{\infty }{{{n}_{i}}}, \quad\left( n \right)!=\prod\limits_{i=0}^{\infty }{{{n}_{i}}}!,  \quad n=\sum\limits_{i=0}^{\infty }{{{n}_{i}}}{{q}^{i}}, \quad 0\le {{n}_{i}}<q,$$
$$\left( n \right)!\left[ {{x}^{n}} \right]{}^{q}{{\varepsilon }^{\left( \varphi  \right)}}\left( x \right)\circ {}^{q}{{\varepsilon }^{\left( \beta  \right)}}\left( x \right)={{\left( \varphi +\beta  \right)}^{\left\{ n \right\}}}=\sum\limits_{m=0}^{n}{{{\left( \begin{matrix}
   n  \\
   m  \\
\end{matrix} \right)}_{q}}}{{\varphi }^{\left\{ m \right\}}}{{\beta }^{\left\{ n-m \right\}}},$$
$${{\left( \begin{matrix}
   n  \\
   m  \\
\end{matrix} \right)}_{q}}=\left\{ \begin{matrix}
   \frac{\left( n \right)!}{\left( m \right)!\left( \left( n-m \right) \right)!},n\left( \bmod {{q}^{k}} \right)\ge m\left( \bmod {{q}^{k}} \right),  \\
   0,n\left( \bmod {{q}^{k}} \right)<m\left( \bmod {{q}^{k}} \right).  \\
\end{matrix} \right.$$      
Let ${}_{\left( 1 \right)}^{q}\varepsilon \left( x \right)$ denote the series associated with $^{q}\varepsilon \left( x \right)$ according to Theorem 5.2.
Then
 $$\left( n \right)!\left[ {{x}^{n}} \right]{}_{\left( 1 \right)}^{q}{{\varepsilon }^{\left( \varphi  \right)}}\left( x \right)=\varphi {{\left( \varphi +n \right)}^{\left\{ n \right\}-1}}.$$
The identity ${}_{\left( 1 \right)}^{q}{{\varepsilon }^{\left( \varphi +\beta  \right)}}\left( x \right)={}_{\left( 1 \right)}^{q}{{\varepsilon }^{\left( \varphi  \right)}}\left( x \right)\circ {}_{\left( 1 \right)}^{q}{{\varepsilon }^{\left( \beta  \right)}}\left( x \right)$  gives an analog of the generalized binomial Abel formula:
$$\left( \varphi +\beta  \right){{\left( \varphi +\beta +n \right)}^{\left\{ n \right\}-1}}=\sum\limits_{m=0}^{n}{{{\left( \begin{matrix}
   n  \\
   m  \\
\end{matrix} \right)}_{q}}\varphi {{\left( \varphi +m \right)}^{\left\{ m \right\}-1}}\beta {{\left( \beta +n-m \right)}^{\left\{ n-m \right\}-1}}}.$$
Since
$$\left[ {{x}^{n}} \right]\left( 1-{{\left( \log \circ {}_{\left( 1 \right)}^{q}\varepsilon \left( x \right) \right)}^{\prime }} \right)\circ {}_{\left( 1 \right)}^{q}{{\varepsilon }^{\left( \varphi  \right)}}\left( x \right)=\left[ {{x}^{n}} \right]{}^{q}{{\varepsilon }^{\left( \varphi +n \right)}}\left( x \right),$$
it follows from
$$\left( 1-{{\left( \log \circ {}_{\left( 1 \right)}^{q}\varepsilon \left( x \right) \right)}^{\prime }} \right)\circ {}_{\left( 1 \right)}^{q}{{\varepsilon }^{\left( \varphi +\beta  \right)}}\left( x \right)=\left( 1-{{\left( \log \circ {}_{\left( 1 \right)}^{q}\varepsilon \left( x \right) \right)}^{\prime }} \right)\circ {}_{\left( 1 \right)}^{q}{{\varepsilon }^{\left( \varphi  \right)}}\left( x \right)\circ {}_{\left( 1 \right)}^{q}{{\varepsilon }^{\left( \beta  \right)}}\left( x \right)$$
that
$${{\left( \varphi +\beta +n \right)}^{\left\{ n \right\}}}=\sum\limits_{m=0}^{n}{{{\left( \begin{matrix}
   n  \\
   m  \\
\end{matrix} \right)}_{q}}{{\left( \varphi +m \right)}^{\left\{ m \right\}}}}\beta {{\left( \beta +n-m \right)}^{\left\{ n-m \right\}-1}}.$$
The identities
$${}_{\left( 1 \right)}^{q}{{\varepsilon }^{\left( \varphi  \right)}}\left( x \right)={{\left( 1,{}_{\left( 1 \right)}^{q}\varepsilon \left( x \right) \right)}_{0}}{}^{q}{{\varepsilon }^{\left( \varphi  \right)}}\left( x \right),  \qquad{}^{q}{{\varepsilon }^{\left( \varphi  \right)}}\left( x \right)={{\left( 1,{}^{q}{{\varepsilon }^{\left( -1 \right)}}\left( x \right) \right)}_{0}}{}_{\left( 1 \right)}^{q}{{\varepsilon }^{\left( \varphi  \right)}}\left( x \right)$$
gives analogs of other Abel identities [18, pp. 92-99], [21]:
$$\varphi {{\left( \varphi +n \right)}^{\left\{ n \right\}-1}}=\sum\limits_{m=0}^{n}{{{\left( \begin{matrix}
   n  \\
   m  \\
\end{matrix} \right)}_{q}}}{{\varphi }^{\left\{ m \right\}}}m{{n}^{\left\{ n-m \right\}-1}},$$
$${{\varphi }^{\left\{ n \right\}}}=\sum\limits_{m=0}^{n}{{{\left( \begin{matrix}
   n  \\
   m  \\
\end{matrix} \right)}_{q}}\varphi {{\left( \varphi +m \right)}^{\left\{ m \right\}-1}}}{{\left( -m \right)}^{\left\{ n-m \right\}}}.$$
In particular, if $q=2$, $n={{2}^{k}}-1$, then
$${{\left( \varphi +{{2}^{k}}-1 \right)}^{k-1}}=\sum\limits_{m=0}^{{{2}^{k}}-1}{{{\varphi }^{\left\{ m \right\}-1}}m{{\left( {{2}^{k}}-1 \right)}^{\left\{ {{2}^{k}}-1-m \right\}-1}}},$$
$${{\varphi }^{k-1}}=\sum\limits_{m=0}^{{{2}^{k}}-1}{{{\left( \varphi +m \right)}^{\left\{ m \right\}-1}}}{{\left( -m \right)}^{\left\{ {{2}^{k}}-1-m \right\}}}.$$
The coefficients ${n\choose m}_{q}$ are elements of the matrix ${{T}^{\left( q \right)}}={{P}_{\left[ 0,q \right]}}\times {{P}_{{}^{q}\varepsilon \left( x \right)}}$ considered in [12]. Note that since
$^{q}\varepsilon \left( x \right)\circ {}^{q}\varepsilon \left( x \right)={}^{q}{{\varepsilon }^{\left( 2 \right)}}\left( x \right)$,  ${{x}^{q}}{\varepsilon }'\left( x \right)\circ {}^{q}{{\varepsilon }^{\left( -1 \right)}}\left( x \right)=x{{\left( \log \circ {}^{q}\varepsilon \left( x \right) \right)}^{\prime }}$,
then
$$\sum\limits_{m=0}^{n}{{{\left( \begin{matrix}
   n  \\
   m  \\
\end{matrix} \right)}_{q}}}={{2}^{\left\{ n \right\}}},   \qquad\sum\limits_{m=0}^{n}{{{\left( \begin{matrix}
   n  \\
   m  \\
\end{matrix} \right)}_{q}}m{{\left( -1 \right)}^{\left\{ n-m \right\}}}}=\left\{ \begin{matrix}
   {{q}^{k}},n={{q}^{k}},  \\
   0,n\ne {{q}^{k}}.  \\
\end{matrix} \right.$$

In the Riordan matrices theory, the concept of pseudo-involution in the Riordan group occupies an important place [5,6,9,17]. Let us introduce a similar concept for the group $R\left( {{P}_{0}} \right)$.
The matrix ${{\left( b\left( x \right),a\left( x \right) \right)}_{0}}$ with the property 
$$\left( b\left( x \right),a\left( x \right) \right)_{0}^{-1}={{\left( 1,-1 \right)}_{0}}{{\left( b\left( x \right),a\left( x \right) \right)}_{0}}{{\left( 1,-1 \right)}_{0}}={{\left( b\left( -x \right),a\left( -x \right) \right)}_{0}},$$
will be called a pseudo-involution in the group  $R\left( {{P}_{0}} \right)$. 
Obviously, if the matrix ${{\left( b\left( x \right),a\left( x \right) \right)}_{0}}$ is a pseudo-involution, then the matrices 
 ${{\left( b\left( -x \right),-a\left( -x \right) \right)}_{0}}$, ${{\left( b\left( x \right),-a\left( x \right) \right)}_{0}}$ are involutions.\\ 
{\bfseries Theorem 5.3. } \emph{If the matrix  ${{\left( b\left( x \right),-a\left( x \right) \right)}_{0}}$ is an involution, then it can be represented as}
$${{\left( b\left( x \right),-a\left( x \right) \right)}_{0}}={{\left( {{b}^{\left( {1}/{2}\; \right)}}\left( x \right),{{a}^{\left( {1}/{2}\; \right)}}\left( x \right) \right)}_{0}}{{\left( 1,-1 \right)}_{0}}\left( {{b}^{\left( {1}/{2}\; \right)}}\left( x \right),{{a}^{\left( {1}/{2}\; \right)}}\left( x \right) \right)_{0}^{-1}.$$
{\bfseries Proof.} If
$$\left( 1,{{a}^{\left( {1}/{2}\; \right)}}\left( x \right) \right)_{0}^{-1}={{\left( 1,{{h}^{\left( -1 \right)}}\left( x \right) \right)}_{0}}, \qquad\left( 1,h\left( x \right) \right)_{0}^{-1}={{\left( 1,{{c}^{\left( {1}/{2}\; \right)}}\left( x \right) \right)}_{0}},$$
then
$${{\left( 1,{{a}^{\left( {1}/{2}\; \right)}}\left( x \right) \right)}_{0}}{{\left( 1,h\left( x \right) \right)}_{0}}={{\left( 1,a\left( x \right) \right)}_{0}},$$ 
$${{\left( 1,{{c}^{\left( {1}/{2}\; \right)}}\left( x \right) \right)}_{0}}{{\left( 1,{{h}^{\left( -1 \right)}}\left( x \right) \right)}_{0}}={{\left( 1,c\left( x \right) \right)}_{0}}, \qquad\left( 1,a\left( x \right) \right)_{0}^{-1}={{\left( 1,c\left( x \right) \right)}_{0}}.$$
Since the matrix ${{\left( 1,a\left( x \right) \right)}_{0}}$  is a pseudo-involution, it follows from the condition $c\left( x \right)=a\left( -x \right)$ that ${{h}^{\left( -1 \right)}}\left( x \right)=h\left( -x \right)$. Then
$${{\left( 1,-a\left( x \right) \right)}_{0}}={{\left( 1,{{a}^{\left( {1}/{2}\; \right)}}\left( x \right) \right)}_{0}}{{\left( 1,-1 \right)}_{0}}{{\left( 1,{{h}^{\left( -1 \right)}}\left( x \right) \right)}_{0}},$$ 
$${{\left( b\left( x \right),-a\left( x \right) \right)}_{0}}={{\left( {{b}^{\left( {1}/{2}\; \right)}}\left( x \right),1 \right)}_{0}}{{\left( 1,-a\left( x \right) \right)}_{0}}{{\left( {{b}^{\left( -{1}/{2}\; \right)}}\left( x \right),1 \right)}_{0}}.  \qquad\square $$

The following theorem concerns unipotents ${{w}_{q}}\left( x \right)\in \left[\!\left[ {{P}_{0,q}},a\left( x \right) \right]\!\right]$, ${{w}_{q}}\left( x \right)=1+{{\eta }_{q}}\left( x \right)$, where the series  ${{\eta }_{q}}\left( x \right)$ is defined by formula (2).\\
{\bfseries Theorem 5.4.} \emph{ The matrices ${{\left( {{w}_{q,i}}\left( x \right),{{w}_{q,j}}\left( x \right) \right)}_{0}}$ form a commutative subgroup in $R\left( {{P}_{0,q}} \right)$ whose elements are multiplied by the rule }
$${{\left( {{w}_{q,1}}\left( x \right),{{w}_{q,2}}\left( x \right) \right)}_{0}}{{\left( {{w}_{q,3}}\left( x \right),{{w}_{q,4}}\left( x \right) \right)}_{0}}={{\left( {{w}_{q,1}}\left( x \right)\circ {{w}_{q,3}}\left( x \right),{{w}_{q,2}}\left( x \right)\circ {{w}_{q,4}}\left( x \right) \right)}_{0}}.$$
{\bfseries Proof.}  If $\left[ {{x}^{n}} \right]{{w}_{q,i}}\left( x \right)\ne 0$, $n>0$, then 
$${{x}^{n}}\circ {{w}_{q,i}}\left( x \right)={{x}^{n}}\circ w_{q,i}^{\left( n \right)}\left( x \right)={{x}^{n}},  \qquad{{w}_{q,i}}{{\left( {{w}_{q,j}}\left( x \right) \right)}_{0}}={{w}_{q,i}}\left( x \right),$$
$${{\left( {{w}_{q,1}}\left( x \right),{{w}_{q,2}}\left( x \right) \right)}_{0}}{{\left( {{w}_{q,3}}\left( x \right),{{w}_{q,4}}\left( x \right) \right)}_{0}}
={{\left( 1+{{\eta }_{q,1}}\left( x \right)+{{\eta }_{q,3}}\left( x \right),1+{{\eta }_{q,2}}\left( x \right)+{{\eta }_{q,4}}\left( x \right) \right)}_{0}}.  \quad\square $$ 

We denote this subgroup by $U\left( {{P}_{0,q}} \right)$.\\
{\bfseries Remark 5.1.} For any group $R\left( {{P}_{0}}\times {{P}_{0,q}} \right)$, the product of the matrices $\left( {{w}_{q,i}}\left( x \right)|{{P}_{0}}\times {{P}_{0,q}} \right)$ is determined by their zero columns, so if
$$\left( {{w}_{q,1}}\left( x \right)|{{P}_{0,q}} \right)\left( {{w}_{q,2}}\left( x \right)|{{P}_{0,q}} \right)=\left( {{w}_{q,3}}\left( x \right)|{{P}_{0,q}} \right),$$
then
$$\left( {{w}_{q,1}}\left( x \right)|{{P}_{0}}\times {{P}_{0,q}} \right)\left( {{w}_{q,2}}\left( x \right)|{{P}_{0}}\times {{P}_{0,q}} \right)=\left( {{w}_{q,3}}\left( x \right)|{{P}_{0}}\times {{P}_{0,q}} \right).$$
Thus, the group $R\left( {{P}_{0}}\times {{P}_{0,q}} \right)$ contains a subgroup isomorphic to the subgroup $U\left( {{P}_{0,q}} \right)$.\\
{\bfseries Theorem 5.5.} \emph{Pseudo-involutions in the group $R\left( {{P}_{0}}\times {{P}_{0,2}} \right)$ form a subgroup isomorphic to the subgroup  $U\left( {{P}_{0,2}} \right)$.}\\
{\bfseries Proof.} If the matrix ${{\left( b\left( x \right),a\left( x \right) \right)}_{0}}$ is a pseudo-involution (excluding matrices ${{\left( 1,-1 \right)}_{0}}$, ${{\left( -1,-1 \right)}_{0}}$, which are both involutions and pseudo-involutions), then
$${{\left( 1,a\left( x \right) \right)}_{0}}b\left( -x \right)={{b}^{\left( -1 \right)}}\left( x \right), \qquad{{\left( 1,a\left( x \right) \right)}_{0}}={{\left( 1,{{a}^{\left( {1}/{2}\; \right)}}\left( x \right) \right)}_{0}}{{\left( 1,h\left( x \right) \right)}_{0}},$$
$$\left( 1,{{a}^{\left( {1}/{2}\; \right)}}\left( x \right) \right)_{0}^{-1}={{\left( 1,{{h}^{\left( -1 \right)}}\left( x \right) \right)}_{0}},  \qquad{{h}^{\left( -1 \right)}}\left( x \right)=h\left( -x \right).$$
Let ${{P}_{0}}={{P}_{0,2}}$.
Since, in the algebra  $\left[\!\left[ {{P}_{0,2}},a\left( x \right) \right]\!\right]$, the series $\eta \left( x \right)=\sum\nolimits_{n=0}^{\infty }{{{\eta }_{2n+1}}}{{x}^{2n+1}}$ is a nilpotent of degree 2, then the series with the property  ${{h}^{\left( -1 \right)}}\left( x \right)=h\left( -x \right)$ is a unipotent. Hence, if the matrix ${{\left( 1,a\left( x \right) \right)}_{0}}$ is a pseudo-involution, then the series $a\left( x \right)$ is a unipotet and
${{\left( 1,a\left( x \right) \right)}_{0}}{{x}^{2n+1}}={{x}^{2n+1}}$. Denote $b\left( x \right)={{b}_{1}}\left( {{x}^{2}} \right)+x{{b}_{2}}\left( {{x}^{2}} \right)$. Then ${{\left( 1,a\left( x \right) \right)}_{0}}b\left( -x \right)={{\tilde{b}}_{1}}\left( {{x}^{2}} \right)-x{{b}_{2}}\left( {{x}^{2}} \right)$. The equation
$${{b}_{1}}\left( {{x}^{2}} \right)\circ {{\tilde{b}}_{1}}\left( {{x}^{2}} \right)+{{b}_{1}}\left( {{x}^{2}} \right)\circ x{{b}_{2}}\left( {{x}^{2}} \right)-{{\tilde{b}}_{1}}\left( {{x}^{2}} \right)\circ x{{b}_{2}}\left( {{x}^{2}} \right)=1$$
has the only solution ${{b}_{1}}\left( {{x}^{2}} \right)={{\tilde{b}}_{1}}\left( {{x}^{2}} \right)=1$. Thus, the series $b\left( x \right)$ is also a unipotet and, hence, ${{\left( b\left( x \right),a\left( x \right) \right)}_{0}}\in U\left( {{P}_{0,2}} \right)$. On the other hand, each element of the group $U\left( {{P}_{0,2}} \right)$ is a pseudo-involution: $\left( {{w}_{q,i}}\left( x \right),{{w}_{q,j}}\left( x \right) \right)_{0}^{-1}={{\left( {{w}_{q,i}}\left( -x \right),{{w}_{q,j}}\left( -x \right) \right)}_{0}}$. Thus, the subgroup of the group $R\left( {{P}_{0}}\times {{P}_{0,2}} \right)$ isomorphic to the subgroup $U\left( {{P}_{0,2}} \right)$ consists of pseudo-involutions of the group $R\left( {{P}_{0}}\times {{P}_{0,2}} \right)$.       \qquad   $\square $

Examples of the group $R\left( {{P}_{0}}\times {{P}_{0,2}} \right)$ are the groups $R\left( {{P}_{\left[ 0,2 \right]}} \right)$ and $R\left( {{P}_{g\left( -1,x \right)}} \right)$.\\
{\bfseries Example 5.3.} Let $_{\left( 1 \right)}h\left( x \right)$ denote the series associated with $h\left( x \right)$ according to  Theorem 5.2. Then if ${{h}^{\left( -1 \right)}}\left( x \right)=h\left( -x \right)$, then matrices $\left( 1,{}_{\left( 1 \right)}{{h}^{\left( 2 \right)}}\left( x \right) \right)_{0}$, $\left( {}_{\left( 1 \right)}{{h}^{\left( 2 \right)}}\left( x \right),{}_{\left( 1 \right)}{{h}^{\left( 2 \right)}}\left( x \right) \right)_{0}$, are pseudo-involutions. Based on this, we construct a pseudo-involution in the group $R\left( {{P}_{0,3}} \right)$. Let $h\left( x \right)={{\left[ {{e}^{x}} \right]}_{3}}{{e}^{{{x}^{3}}}}=\left( 1+x+{{{x}^{2}}}/{2}\; \right){{e}^{{{x}^{3}}}}$. Then, taking into account identity (3),
$${{h}^{\left( \varphi  \right)}}\left( x \right)=\sum\limits_{n=0}^{\infty }{\frac{{{\varphi }^{n}}}{n!}}{{x}^{3n}}+\sum\limits_{n=0}^{\infty }{\frac{{{\varphi }^{n+1}}}{n!}}{{x}^{3n+1}}+\sum\limits_{n=0}^{\infty }{\frac{{{\varphi }^{n+2}}}{n!2}}{{x}^{3n+2}},$$
$$_{\left( 1 \right)}{{h}^{\left( \varphi  \right)}}\left( x \right)=\sum\limits_{n=0}^{\infty }{\frac{\varphi {{\left( \varphi +3n \right)}^{n-1}}}{n!}}{{x}^{3n}}+\sum\limits_{n=0}^{\infty }{\frac{\varphi {{\left( \varphi +3n+1 \right)}^{n}}}{n!}}{{x}^{3n+1}}
+\sum\limits_{n=0}^{\infty }{\frac{\varphi {{\left( \varphi +3n+2 \right)}^{n+1}}}{n!2}}{{x}^{3n+2}},$$
$$\left( {}_{\left( 1 \right)}{{h}^{\left( 2 \right)}}\left( x \right),{}_{\left( 1 \right)}{{h}^{\left( 2 \right)}}\left( x \right) \right)_{0}
=\left( \begin{matrix}
   1 & 0 & 0 & 0 & 0 & 0 & 0 & 0 & 0 & \cdots   \\
   2 & 1 & 0 & 0 & 0 & 0 & 0 & 0 & 0 & \cdots   \\
   4 & 4 & 1 & 0 & 0 & 0 & 0 & 0 & 0 & \cdots   \\
   2 & 0 & 0 & 1 & 0 & 0 & 0 & 0 & 0 & \cdots   \\
   12 & 4 & 0 & 8 & 1 & 0 & 0 & 0 & 0 & \cdots   \\
   49 & 32 & 6 & 40 & 10 & 1 & 0 & 0 & 0 & \cdots   \\
   8 & 0 & 0 & 8 & 0 & 0 & 1 & 0 & 0 & \cdots   \\
   81 & 20 & 0 & 96 & 10 & 0 & 14 & 1 & 0 & \cdots   \\
   500 & 242 & 36 & 676 & 140 & 12 & 112 & 16 & 1 & \cdots   \\
   \vdots  & \vdots  & \vdots  & \vdots  & \vdots  & \vdots  & \vdots  & \vdots  & \vdots  & \ddots   \\
\end{matrix} \right).$$

E-mail: {evgeniy\symbol{"5F}burlachenko@list.ru}
\end{document}